\documentclass[a4paper,11pt]{amsart}

\usepackage{amsfonts,latexsym,rawfonts,amsmath,amssymb,amsthm, mathrsfs, lscape}
\usepackage{stackrel}
\usepackage[all]{xy}
\usepackage[english]{babel}
\usepackage[utf8]{inputenc}
\usepackage{xfrac}
\usepackage{multirow}
\usepackage{enumitem}
\usepackage[usenames, dvipsnames]{xcolor}
\RequirePackage{tikz, tikz-cd}

\usepackage{array, tabularx}

\usepackage{setspace}

\usepackage{textcomp}

\usepackage[hypertexnames=false,
    pdftex,
    pdfpagemode=UseNone,
    breaklinks=true,
    extension=pdf,
    colorlinks=true,
    linkcolor=blue,
    citecolor=blue,
    urlcolor=blue,
]{hyperref}
\usepackage[top=35mm, bottom=35mm, left=3.5cm, right=3.5cm, headsep=6mm, footskip=10mm, headheight=2cm]{geometry}

\renewcommand{\Re}{\mathsf{Re}\,}
\renewcommand{\Im}{\mathsf{Im}\,}

\newcommand{\rk}{\mathrm{rank}}

\newcommand{\Pic}{\mathrm{Pic}}
\newcommand{\ODA}{\overline{\DA\vphantom{A} }}
\newcommand{\Def}{\mathrm{Def}}
\newcommand{\Lg}{\mathfrak{g}}

\newcommand{\ov}{\overline}

\newcommand{\e}{\mathrm{e}}

\newcommand{\del}{\partial}

\newcommand{\Aut}{\mathrm{Aut}}

\newcommand{\id}{\operatorname{\text{\sf id}}}
\renewcommand{\Re}{\operatorname{Re}}
\renewcommand{\Im}{\operatorname{Im}}
\newcommand{\im}{\operatorname{im}}
\newcommand{\kod}{\operatorname{kod}}

\newcommand{\Oo}{\mathcal{O}}

\newcommand{\la}{\lambda}

\newcommand{\ce}{\mathcal{C}^\infty}
\newcommand{\CC}{\mathbb{C}}

\newcommand{\RR}{\mathbb{R}}

\newcommand{\ZZ}{\mathbb{Z}}
\newcommand{\NN}{\mathbb{N}}
\newcommand{\PP}{\mathbb{P}}

\newcommand{\QQ}{\mathbb{Q}}

\newcommand{\Ka}{K\"{a}hler}

\newcommand*\DA{\mathop{}\!\mathbin\Box}

\def\XXint#1#2#3{{\setbox0=\hbox{$#1{#2#3}{\int}$ }
\vcenter{\hbox{$#2#3$ }}\kern-.6\wd0}}

\allowdisplaybreaks[1]

\newtheorem{theorem}{Theorem}[section]
\newtheorem*{theorem*}{Theorem}
\newtheorem{theoremA}{Theorem}
\newtheorem{propA}{Proposition}

\addtocounter{theoremA}{+1}

\newtheorem{corollary}[theorem]{Corollary}
\newtheorem*{corollary*}{Corollary}
\newtheorem{lemma}[theorem]{Lemma}
\newtheorem{definition}[theorem]{Definition}
\newtheorem{proposition}[theorem]{Proposition}
\newtheorem{remark}[theorem]{Remark}

\newtheorem{claim}{Claim}

\newtheoremstyle{mytheoremstyle} 
   {}
    {}
    {}
    {}                           
    {\bfseries}
    {.}                          
    {.5em}                       
    {}  

\theoremstyle{mytheoremstyle}
\newtheorem{example}{Example}[section]

\title[Vaisman manifolds with $c_1=0$]{Vaisman manifolds with vanishing first Chern class}

\author{Nicolina Istrati}
\address{Nicolina Istrati\\FB 12/Mathematik und Informatik\\
Philipps-Universit\"at Marburg\\
Hans-Meer\-wein-Str. 6\\
35032 Marburg\\
Germany}
\email{istrati@mathematik.uni-marburg.de}

\begin{document}

\begin{abstract}
Compact Vaisman manifolds with vanishing first Chern class split into three categories, depending on the sign of the Bott-Chern class. We show that Vaisman manifolds with non-positive Bott-Chern class admit canonical metrics, are quasi-regular and are stable under deformations. We also show that Calabi-Yau Vaisman manifolds satisfy a version of the Beauville-Bogomolov decomposition and have torsion canonical bundle.  Finally, we prove a general result concerning the behaviour of the automorphism group of a complex manifold under deformations.
\end{abstract}
\maketitle

\section{Introduction}

Kähler manifolds with vanishing first Chern class have remarkable properties. They have canonical metrics \cite{Yau}, i.e. unique Ricci-flat Kähler metrics with prescribed cohomology class. They split, after a finite cover, into a product of a torus and a simply connected manifold with trivial canonical bundle \cite{Be84}, so in particular they have torsion canonical class. Finally, they are stable under small deformations and have smooth, universal Kuranishi spaces.

All of these properties fail as soon as one leaves the Kähler realm. First of all, as advocated by Tosatti \cite{To}, the reasonable condition to ask in this generality is the vanishing of the first Chern class in Bott-Chern cohomology, which refines the de Rham cohomology. We will call such manifolds (non-Kähler) Calabi-Yau. In this case, one can still show the existence of  certain distinguished metrics given as solutions of Monge-Ampère equations \cite{TW}, however their unicity can no longer be prescribed by topological data. This phenomenon is related to the fact that the identity component of the automorphism group of such manifolds no longer is compact, and therefore cannot be an isometry group - this happens for instance for all complex parallelizable nilmanifolds.

Non-Kähler Calabi-Yau manifolds also fail to have holomorphically torsion canonical bundle, as illustrated by certain deformations of complex parallelizable solvmanifolds \cite[Ex. III-(3b)]{Na}, \cite[Ex.~3.2]{To}. In particular, one should not expect any reasonable Beauville-Bogomolov decomposition, even for $\del\ov\del$-manifolds with trivial canonical bundle, as was remarked in \cite[Ex. 2.14]{ACRT}. Finally, manifolds with trivial canonical bundle can have obstructed Kuranishi spaces \cite{Ro}.

On the other hand, Vaisman manifolds \cite{V'} are a class of Hermitian manifolds which remain rather close to the Kähler world. Indeed, while they no longer satisfy the $\del\ov\del$-lemma, they still have a good cohomological behaviour. More importantly, they are endowed with a holomorphic one-dimensional transversally Kähler foliation. One of the main goals of the present paper is to show  that the Vaisman Calabi-Yau manifolds have properties analogous to the Kähler ones.

Let $X=(M,J)$ be a compact complex manifold with $\dim_\CC X>1$. A Hermitian metric on $X$ with fundamental form $\Omega$ is called Vaisman if the following conditions hold:
\[d\Omega=\theta\wedge\Omega,\quad \nabla\theta=0\]
where $\theta$ is the \textit{Lee form}, uniquely determined by the metric, 
and $\nabla$ is the Levi-Civita connection. The vector field which is metric dual to $\theta$ is real holomorphic and Killing \cite{V}, thus generating a complex one-dimensional subgroup $G\subset \Aut(X)$, which 
 is independent of the chosen Vaisman metric  \cite{Ts97}. $X$ is called quasi-regular if $G$ is compact, in which case the complex space $X/G$ has a natural structure of a complex projective orbifold.

The starting point of our discussion is the remark that compact Vaisman manifolds with vanishing first Chern class in de Rham cohomology split into three classes. Indeed, by \cite{IO22}, the kernel of the natural map 
\begin{equation}\label{eq: BCdR}H^{1,1}_{BC}(X,\RR)\rightarrow H^2(X,\RR)\end{equation}
is one-dimensional, generated by the class of a transverse Kähler metric. Consequently, a Vaisman manifold with $c_1(X)=0$ in $H^2(X,\RR)$ has a well defined sign, depending on whether $c_1^{BC}(X)=0$ or $\pm c_1^{BC}(X)$ can be represented by a transverse Kähler class, which we denote by $\pm c_1^{BC}(X)>0$. The three classes can be distinguished by their Kodaira dimension  $\kod(X)$:

\begin{propA}
Let $X$ be an $n$-dimensional  compact Vaisman manifold with $c_1(X)=0\in H^2(X,\RR)$. 
\begin{enumerate}
\item If $c_1^{BC}>0$ then $\kod(X)=-\infty$. 
\item If $c_1^{BC}=0$ then $\kod(X)=0$. 
\item If $c_1^{BC}<0$ then $\kod(X)=n-1$.
\end{enumerate}
\end{propA}

We focus our study on the two classes $c_1^{BC}(X)=0$ and $c_1^{BC}(X)<0$, which behave quite similarly. First of all, we obtain analogous versions of the Aubin-Yau theorem, producing canonical Vaisman metrics with prescribed topological data. Note that in the non-Kähler setting, the role of the Ricci curvature of a metric $\Omega$  is taken by the Chern-Ricci form  $\rho_\Omega:=-i\del\ov\del\log\det\Omega$  which represents  $2\pi c_1^{BC}(X)$. 
Moreover, we have a natural map 
\[\psi:H^1(X,\RR)\rightarrow H^{1,1}_{BC}(X,\RR), \quad [\theta]\mapsto[-dJ\theta]_{BC}.\]

\begin{theoremA}\label{thm: A}
Let $X$ be a compact Vaisman manifold with $c_1(X)=0$ and let $\tau\in H^1(X)$ represent some Lee form of a Vaisman metric on $X$.
\begin{enumerate}
\item If $c_1^{BC}(X)=0$, then there exists a unique Chern-Ricci flat Vaisman metric with Lee class $\tau$ and prescribed total volume.
\item If $c_1^{BC}(X)<0$ and $\tau\in\psi^{-1}(-2\pi c_1^{BC}(X))$, then there exists a unique  Vaisman metric $(\Omega,\theta)$ with $\theta\in\tau$, $\rho_\Omega=dJ\theta$ and prescribed total volume.
\end{enumerate}
\end{theoremA}

Remark that in the negative case, the curvature condition  $\rho_\Omega=dJ\theta$ plays the role of the Einstein equation, and that $\psi(\tau)=-2\pi c_1^{BC}(X)$ is a necessary cohomological condition to ensure such prescribed Chern-Ricci curvature. At the same time,  $\psi^{-1}(-2\pi c_1^{BC}(X))$ always contains a Lee class, so canonical metrics do indeed exist. 

Related to the existence of canonical metrics, we also obtain a rigidity result for these manifolds:

\begin{theoremA}\label{thm: B}
Let $X$ be a compact Vaisman manifold with $c_1^{BC}(X)\leq 0$. Then $\Aut_0(X)$ is compact one-dimensional and $X$ is quasi-regular.
\end{theoremA}

For Calabi-Yau Vaisman manifolds, we further obtain an analogue of the Beauville-Bogomolov decomposition (see Theorem~\ref{thm: holoVF} for a more precise statement). Note that Vaisman manifolds cannot split holomorphically into products of complex manifolds \cite[Thm.~5.1]{V}, so the following result is the best one could hope for.

\begin{theoremA}\label{thm: C}
Let $X$ be a compact Vaisman Calabi-Yau manifold. 
Then there exists a finite étale cover of $X$ which is a definite principal  elliptic orbi-bundle over an orbifold $Y$ which is the product of an abelian variety and a simply connected orbifold with trivial canonical divisor. In particular, $X$ has holomorphically torsion canonical bundle.
\end{theoremA}

Regarding the deformation theory, it should be noted that the Vaisman condition is not preserved under deformations, and neither is the quasi-regular condition under Vaisman deformations. Both of these phenomena are generic and can already be detected on Hopf surfaces, which in particular satisfy $c_1^{BC}(X)>0$. Nonetheless, this turns out to be an effect of the presence of too many symmetries, cf. Corollary \ref{cor: defnVaisman}. In particular, the conditions $c_1^{BC}(X)\leq 0$  impose enough robustness, to the effect that:

\begin{theoremA}\label{thm: D}
Any small deformation $X_t$ of  a compact Vaisman manifold $X$ with $c_1^{BC}(X)\leq 0$ is again Vaisman with $c_1^{BC}(X_t)\leq 0$.
\end{theoremA}

This  statement is deduced from a more general result which we prove, which has an independent interest:

\begin{theoremA}\label{thm: E}
Let $X$ be a compact complex manifold and let $H\subset\Aut_0(X)$ be a compact complex subgroup. Then for any small deformation $X_t$ of $X$ one has $\dim_\CC\Aut(X_t)\geq\dim_\CC H$. In particular, if  $\Aut_0(X)$ is compact, then for any small deformation $X_t$ of $X$, $\Aut_0(X_t)$ is a small deformation of $\Aut_0(X)$. 
\end{theoremA}

This is shown by first noting that the group $H$ must act trivially on the relevant cohomology groups. Secondly, one proves that $H$ acts trivially also on the Kuranishi space, using the description of the latter via Hodge theory.

Theorems \ref{thm: A} and \ref{thm: B}  are proven in different manners for the null and the negative case. The existence of canonical metrics on Vaisman Calabi-Yau manifolds is deduced from the the Yau theorem for Vaisman metrics of Ornea-Verbitsky \cite[Thm~4.1]{ov22}, which is a  based on the transverse Yau theorem of El Kacimi-Alaoui \cite[Sect.~3.5.5]{ka}. As a consequence of the uniqueness of the canonical metrics, we immediately obtain that the group $\Aut_0(X)$ of such manifolds is compact, and then deduce, via the description of the automorphism groups of locally conformally Kähler manifolds  \cite[Cor.~1]{ist19}, that it is one-dimensional. 

In the negative case, we first prove that such manifolds are quasi-regular. This is done via Weitzenböck techniques involving both the Weyl and the Chern connections, and using again the transverse Yau theorem \cite[Sect.~3.5.5]{ka}. Then we make use of the orbifold version of the Aubin-Yau theorem   \cite[Thm~1.3]{F} to deduce the existence of canonical metrics. 

In particular, Theorem~\ref{thm: B} leads us to the realm of complex orbifolds. In Theorem~\ref{thm: equiv qregV}, we make precise the relation between quasi-regular Vaisman manifolds and  projective orbifolds. Vaguely speaking, the former correspond to positive elliptic orbi-bundles over the latter, and the sign of $c_1^{BC}$  is preserved under this correspondence. 

However, the orbifolds $X/G$ are not generally \textit{pure} in the sense of  \cite{camp}, i.e. the orbifold structure is not determined by the complex structure  of $X/G$ alone, but also by a $\QQ$-divisor $D$, related to the codimension $1$ subspace of $X$ where $G$ acts non-freely. Nonetheless, in the Calabi-Yau setting we can show that, after passing to a finite cover, one does indeed obtain a pure orbifold as a quotient. This builds on the orbifold version of the Cheeger-Gromoll theorem \cite[Thm.~2]{BZ} and a result of Kollar \cite[Prop.~10.2]{ko05}. In particular, we can then use the Beauville-Bogomolov decomposition for pure Calabi-Yau orbifolds \cite[Thm.~6.4]{camp} in order to obtain Theorem~\ref{thm: C}

The Fano case in Kähler geometry is famously difficult, therefore we do not deal with the case $c_1^{BC}>0$ in the present paper. Suffice it to say that none of the above results holds in this setting. As a side remark, the Einstein-Weyl condition \cite{Gau}, which can be seen as a conformal metric analogue of the Calabi-Yau condition, implies $c_1^{BC}>0$. 

It is worth studying the conditions $c_1(X)=0$ or $c_1^{BC}(X)=0$   more generally for locally conformally Kähler manifolds. However, one of the issues is that we lack yet a good description of the cohomology of locally conformally Kähler manifolds, so that in particular the kernel of the map \eqref{eq: BCdR} is unkown in this more general setting. Moreover, we are not aware of any non-Vaisman Calabi-Yau example. If any should exist, they would be very interesting if only from the point of view of locally conformally Kähler geometry.

The paper is organised as follows. In Section~\ref{sec: Orb}, we recall  some of the main facts about complex orbifolds following \cite{BG} and establish the relation between these and quasi-regular Vaisman manifolds. Section~\ref{sec: defn} deals with deformations of complex manifolds and their automorphism groups.  In Section~\ref{sec: c1=0}, we introduce the $c_1=0$ condition on Vaisman manifolds and explain the splitting of such manifolds into the  three classes mentioned before. In Section~\ref{sec: CY}, we prove all the above results for Calabi-Yau Vaisman manifolds. In particular, we also obtain in Section~\ref{sec: inv} a description of the fundamental group of such manifolds, as well as obstructions on the first Betti number.  Section~\ref{sec: negative} deals with the negative case. Finally, in Section~\ref{sec: examples} we present a few explicit examples of Vaisman manifolds with $c_1=0$.

\subsection*{Convention/Notation}
For a given LCK metric $g$, $\Omega$ denotes its fundamental form, $\theta$ its Lee form, $B$ is the metric dual of $\theta$, i.e. the Lee vector field, and we denote by $\omega_0:=-dJ\theta$. For a Vaisman manifold $X$,  $G\subset \Aut(X)$ is the complex Lie group generated by the Lee vector field. 
$\Delta\subset\CC$ will denote the unit disk.

\subsection*{Acknowledgements}
This project started through stimulating conversations with Alexandra Otiman, to which I am indebted. I am also very thankful to Sönke Rollenske for helpful discussions. 

I was partially supported by a grant of Ministry of Research and Innovation, CNCS - UEFISCDI, project no. PN-III-P1-1.1-TE-2021-0228, within PNCDI III, and by  DFG project 509274422.

\section{Complex orbifolds and Vaisman manifolds}\label{sec: Orb}

In this section we recall the main aspects of Vaisman manifolds and of complex orbifolds, while focussing on the relation between the two geometries. For an extensive presentation of locally conformally Kähler and Vaisman geometry, the reader is invited to consult \cite{OV22}, whereas for the theory of orbifolds we follow \cite[Chapter~4]{BG}.

\begin{definition} A Hermitian metric $g$ on a complex manifold $X=(M,J)$ is called \textit{locally conformally Kähler} (LCK) if its fundamental form $\Omega:=g(J\cdot,\cdot)$ satisfies the integrability condition 
    \[d\Omega=\theta\wedge\Omega\] 
    for a closed one-form $\theta$, called \textit{the Lee form} of the metric. If furthermore $\nabla^g\theta=0$, where $\nabla^g$ is the Levi-Civita connection, the metric is called \textit{Vaisman}.
    \end{definition}
    A de Rham class $\tau\in H^1(X,\RR)$ is said to be a \textit{Lee class} if it can be represented by the Lee form of an LCK metric. We will denote by $\mathrm{Lee}(X)\subset H^1(X,\RR)$ the set of all such classes. 
    
  A Vaisman metric is unique in its conformal class up to multiplication by positive constants. Moreover, the norm $||\theta||$ is constant. A Vaisman metric will be called \textit{normalized} if furthermore we impose $||\theta||=1$. In this case, the fundamental form satisfies
    \begin{equation}\label{eq: normV}
    \Omega=-dJ\theta+\theta\wedge J\theta.\end{equation}

When $\dim_\CC X>1$, the Lee form of an LCK metric is uniquely determined by the metric. We denote by $B$ the real vector field which is metric dual to $\theta$, and call it \textit{the Lee vector field}. It is easy to see that if the Lee vector field of an LCK metric is  Killing, then the metric is Vaisman. Conversely, for a Vaisman metric $g$, both $B$ and $JB$ are Killing and real holomorphic \cite{V}, hence they generate a complex one-dimensional subgroup $G\subset \Aut(M,J,g)$, which we will call in the sequel \textit{the Lee group} of the manifold. In fact, by \cite[Cor~2.7]{Ts97}, the group $G$ is independent of the Vaisman metric. 


The group $G$ acts locally freely on $X$, and its orbits define a foliation $\mathcal F$, called \textit{the canonical foliation}. Moreover, any Vaisman metric $(\Omega,\theta)$ has an associated transverse Kähler metric $\omega_0:=-dJ\theta$  with respect to $\mathcal F$ \cite{V}.  If $G$ is compact, then $X$ is called a \textit{quasi-regular Vaisman} manifold. In this case, the quotient $X/G$ inherits a natural structure of a complex orbifold. If moreover $G$ acts freely on $X$, then $X$ is called \textit{a regular Vaisman} manifold, and the quotient $X/G$ is a smooth complex manifold.

Let us recall that a \textit{complex orbifold} $\mathcal Y=(Y,\mathcal U)$ of dimension $n$ is a paracompact Hausdorff space $Y$ endowed with an equivalence class of an orbifold atlas $\mathcal U$, consisting in the following. For each $x\in Y$ there exists a neighbourhood $U\subset Y$ and an orbifold chart $(\tilde U,\Gamma,\phi)\in\mathcal U$ such that $\tilde U\subset\CC^n$ an open set, $\Gamma$ is a finite group of biholomorphisms of $\tilde U$,  and the map $\phi:\tilde U\rightarrow U$ is continuous, $\Gamma$-invariant  and induces a homeomorphism $\tilde U/\Gamma\cong U$. Moreover, the orbifold charts satisfy natural compatibility conditions. Namely, for any point $x\in Y$ covered by two orbifold charts $(\tilde U,\Gamma_U,\phi_U)$ and $(\tilde V,\Gamma_V,\phi_V)$, there exists a third orbifold chart $(\tilde W,\Gamma_W,\phi_W)\in\mathcal U$ with $x\in \phi_W(\tilde W)\subset\phi_U(\tilde U)\cap\phi_V(\tilde V)$ and embeddings $\lambda_{WU},\lambda_{WV}$ making the following diagram commutative 
\begin{equation}\label{eq: orbCh}
\begin{tikzcd}
\tilde U\arrow[d, "\phi_U"] \arrow[r, hookleftarrow, "\lambda_{WU}"] &\tilde W \arrow[r, hook, "\lambda_{WV}"] \arrow[d, "\phi_W"]  &\tilde V\arrow[d,"\phi_V"]\\
\phi_U(\tilde U) \arrow[r, hookleftarrow] &\phi_W(\tilde W) \arrow[r, hook] &\phi_V(\tilde V).
\end{tikzcd}
\end{equation}
Two orbifold atlases $\mathcal U$, $\mathcal V$ are equivalent if they have a common refinement $\mathcal W$, meaning that for any charts $(\tilde U,\Gamma_U,\phi_U)\in\mathcal U$, $(\tilde V,\Gamma_V,\phi_V)\in\mathcal V$ at $x\in Y$, there is a third chart $(\tilde W,\Gamma_W,\phi_W)\in\mathcal W$ at $x$ with embeddings $\lambda_{WU},\lambda_{WV}$ satisfying \eqref{eq: orbCh}. Any orbifold atlas can be refined to a  unique maximal one. We will say that an orbifold structure is  trivial  if it can be represented by an orbifold atlas in which all the groups $\Gamma$ are trivial, i.e. it can be represented by a manifold atlas.


For each $x\in Y$, pick an orbifold chart at $x$ $(\tilde U,\Gamma,\phi)\in\mathcal U$ and an element $\tilde x\in\phi^{-1}(x)$. The conjugacy class of the group $\Gamma_x:=\{g\in \Gamma\mid g\tilde x=\tilde x\}\subset\Gamma$  depends only on $x$ and the equivalence class of $\mathcal U$, and any of its representative is called \textit{the local uniformizing group at $x$}. By the slice theorem, there exists a holomorphic coordinate change sending $\tilde x$ to $0\in\CC^n$ and $\Gamma_x$ to a subgroup of $U(n)$. We will usually tacitly assume that we are working with such linearised orbifold charts.

In particular, the complex orbifold structure induces  a natural structure of a normal complex space on $Y$, by defining the structure sheaf of holomorphic functions $\Oo_Y$ to have stalks $\Oo_{Y,x}=\Oo_{\tilde U,0}^{\Gamma}$, i.e. a local holomorphic function at $x$ is a $\Gamma$-invariant local holomorphic function on $(\tilde U,0)$, for an orbifold chart $\tilde U$ at $x$. However, the underlying analytic structure of $Y$ does not determine the orbifold structure $\mathcal Y$. In fact, as $Y$ is a normal space, its singular locus $Y_{sing}$ has codimension at least $2$, while \textit{the orbifold singular locus} 
\[\Sigma^{orb}(\mathcal Y):=\{x\in Y\mid \Gamma_x\neq\{1\}\}\]
does satisfy $Y_{sing}\subset\Sigma^{orb}(\mathcal Y)$, but the inclusion is strict if some local uniformizing group contains reflections. The codimension-one part of $\Sigma^{orb}(\mathcal Y)$ is encoded in the \textit{branch divisor}:
\begin{equation}\label{eq: branchDiv}
D:=\sum_{\alpha}\left(1-\frac{1}{m_\alpha}\right) D_\alpha
\end{equation}
where the $D_\alpha$'s run over all irreducible Weil divisors contained in $\Sigma^{orb}(\mathcal Y)$, and for each $\alpha$, $m_\alpha:=\mathrm{gcd}_{x\in D_\alpha}|\Gamma_x|$. The $\QQ$-divisor $D$ remembers the branching locus of all the local branched coverings $\phi:\tilde U\rightarrow U$ for $(\tilde U,\Gamma,\phi)\in\mathcal U$. As such, the orbifold structure $\mathcal Y$ is completely determined by the pair $(Y,D)$, where $Y$ is now viewed as a complex space and $D$ is a $\QQ$-divisor on $Y$. If $D=\emptyset$, the orbifold $\mathcal Y$ is also called \textit{pure}, following \cite{camp}.

A tensor on an orbifold $\mathcal Y$ is a tensor $T$ on $Y-\Sigma^{orb}(\mathcal Y)$ so that for each orbifold chart $(\tilde U,\Gamma,\phi)$, $\phi^*T$ extends to a tensor on $\tilde U$. Regarding sheaves and bundles, 
one now has two different notions  on an orbifold. On the one hand, one has the usual sheaf, bundle and divisor theory on the complex space $Y$. As such, a canonical divisor $K_Y$ for instance is defined as any Weil divisor on $Y$ which, restricted to the regular locus $Y_{reg}$, coresponds to the line bundle $\Omega^n_{Y_{reg}}$. 

On the other hand, one also has the richer categories of orbi-sheaves and of orbi-bundles. These are collections of objects living on the uniformizing charts, but which don't generally correspond to any object on $Y$. As such, an \textit{orbi-sheaf}, respectively an \textit{orbi-bundle} $\mathcal E$ on $\mathcal Y$ is given by the collection of $\Gamma$-equivariant sheaves, respectively bundles $\mathcal E_{\tilde U}$ on $\tilde U$, for each $(\tilde U,\Gamma, \phi)\in\mathcal U$, which are compatible with any change of orbifold charts. In particular, the collection $(\Oo_{\tilde U})_{\tilde U\in\mathcal U}$ defines the orbi-sheaf $\Oo_{\mathcal Y}$, the collection $(T\tilde U)_{\tilde U\in\mathcal U}$ defines the tangent orbi-bundle $T\mathcal Y$, while the collection $(K_{\tilde U})_{\tilde U\in\mathcal U}$ defines the canonical orbi-sheaf or orbi-bundle  $K_{\mathcal Y}^{orb}=\det T\mathcal Y^\vee.$ A section of an orbi-sheaf is a compatible collection of $\Gamma$-invariant sections over each orbifold chart.

Note that any divisor $W$ on $Y$ pulls back, via each branched covering map, to an orbi-divisor on $\mathcal Y$, which we denote formally by $\phi^*W$. 
One then has, via the Hurwitz formula \cite[Prop.~4.4.15]{BG}:
\begin{equation}\label{canonical}
    K_{\mathcal Y}^{orb}\equiv\phi^*(K_Y+D)=\phi^* K_Y+R
\end{equation}
where $R$ is the ramification orbi-divisor.

If $H$ is a complex Lie group and $\mathcal P=(P_{\tilde U})_{\tilde U\in\mathcal U}$ is a principal holomorphic $H$ orbi-bundle on $\mathcal Y$, then the total space of $\mathcal P$ has a natural orbifold structure. The underlying complex space $P$ is given by the patching of the charts $\mathcal P_{\tilde U}/\Gamma$, while the orbifold charts are induced by  $(P_{\tilde U},\Gamma,P_{\tilde U}\rightarrow P_{\tilde U}/\Gamma)$. Note that $\Gamma$ acts on the fiber over $0$ of $P_{\tilde U}$, defining a group homomorphism  $h_{\tilde U}:\Gamma\rightarrow H=\Aut(P_{\tilde U,0})$. In particular, $\Gamma$ acts freely on $P_{\tilde U,0}$ precisely when $h_{\tilde U}$ is injective, and in this case the natural orbifold atlas refines to the trivial manifold atlas of $P$.  Finally, note that the natural local projections define a global orbifold fibration map $\pi:\mathcal P\rightarrow \mathcal Y$.

One can then show, via the slice theorem:
\begin{theorem}\label{thm: principal}
If $X$ is any complex manifold and $H$ is a complex Lie group acting holomorphically, properly and locally freely on $X$, then the complex space $Y=X/H$ comes endowed with a natural structure of a complex orbifold $\mathcal Y$ such that $X$ is the underlying complex manifold of an $H$-principal orbi-bundle on $\mathcal Y$. Conversely, if $\mathcal P$ is any holomorphic $H$-principal orbi-bundle on a complex orbifold $\mathcal Y$ such that the defining maps $(h_{\tilde U})_{\tilde U\in\mathcal U}$ are injective, then the total space $P$ is a smooth complex manifold and the orbifold structure of $\mathcal P$ is trivial. 
\end{theorem}

Since orbi-bundles and orbi-sheaves on an orbifold $\mathcal Y=(Y,\mathcal U)$ don't actually correspond to sheaves on $Y$, one cannot directly use the cohomological tools on $Y$ to study them.  Instead, one can define a classifying space $B\mathcal Y$ for $\mathcal Y$ where orbi-sheaves and orbifold cohomology naturally live. Recall that, up to a refinement of $\mathcal U$, all  the groups $\Gamma$ can be seen as subgroups of $U(n)$. Fix a classifying space $BU(n)$ of $U(n)$, together with a universal $U(n)$-bundle $EU(n)\rightarrow BU(n)$. One then defines the CW complex $B\mathcal Y$ as having  local charts
\[\tilde U\times_\Gamma EU(n)\]
for each $(\tilde U,\Gamma,\phi)\in\mathcal U$, where  $\Gamma$ acts on $EU(n)$ via its natural inclusion in $U(n)$. Moreover, the classifying space comes endowed with a natural continuous map 
\[cl:B\mathcal Y\rightarrow Y.\]  
Note that the fiber of the map $cl$ over a point with uniformizing group $\Gamma$ is an Eilenberg-MacLane space $K(\Gamma,1)$.

On $B\mathcal Y$ one can now define \textit{holomorphic} objects in a chart $\tilde U\times_\Gamma EU(n)$ as objects defined on $\tilde U\times EU(n)$ which are $\Gamma$-invariant, holomorphic in the $\tilde U$ variable and smooth in the $EU(n)$-variable, i.e. smooth on any finite model of $EU(n)$. In this way, one shows that holomorphic bundles on $B\mathcal Y$ are in a bijective correspondence with orbi-bundles on $\mathcal Y$ \cite[Thm.~4.3.11]{BG}. In particular, the orbi-sheaf $\Oo_{\mathcal Y}$ corresponds to the sheaf  of holomorphic functions on $B\mathcal Y$, which we still denote by $\Oo_{\mathcal Y}$.

On $B\mathcal Y$ one  has the usual cohomological interpretation. For instance, one defines the orbifold homotopy groups as $\pi_i^{orb}(\mathcal Y):=\pi_i(B\mathcal Y)$. 
For a ring $R$, the orbifold cohomology of $\mathcal Y$ with coefficients in $R$ is defined as
\[H^\bullet_{orb}(\mathcal Y,R):=H^\bullet(B\mathcal Y,R)\]
and if $R$ is a field, the map $cl^*:H^\bullet(Y,R)\rightarrow H^\bullet_{orb}(\mathcal Y,R)$ is an isomorphism. 

If $G=\CC/\Lambda$ is a one-dimensional compact complex torus, consider the sheaf $\Oo_{\mathcal Y}(G)$ of $G$-valued holomorphic functions  on $B\mathcal Y$  and $\Oo^*_{\mathcal Y}\subset\Oo_{\mathcal Y}$ the sheaf of non-vanishing holomorphic functions. Then via the classical theory, one has
\[H^1_{orb}(\mathcal Y,\Oo^*_{\mathcal Y}):=H^1(B\mathcal Y,\Oo^*_{\mathcal Y})\cong \Pic^{orb}(\mathcal Y)\]
where $\Pic^{orb}(\mathcal Y)$ denotes the isomorphism classes of holomorphic line orbi-bundles over $\mathcal Y$. Similarly, 
\[H^1_{orb}(\mathcal Y,\Oo_{\mathcal Y}(G))\cong \{ \text{holomorphic }G\text{-principal orbi-bundle on  }\mathcal Y\}/_{\cong} .\] In particular, the exact sequence of sheaves on $B\mathcal Y$ 
\[    \begin{tikzcd}
0\rar &\Lambda_{B\mathcal Y}\rar &\Oo_{\mathcal Y}\rar &\Oo_{\mathcal Y}(G)\rar &0
\end{tikzcd} \]
defines the integral Chern class of a $G$-principal orbi-bundle as the induced connecting map in cohomology
\[c^\ZZ:H^1_{orb}(\mathcal Y,\Oo_{\mathcal Y}(G))\rightarrow H^2_{orb}(\mathcal Y,\Lambda)\cong H_{orb}^2(\mathcal Y,\ZZ)\otimes\Lambda.\]

\begin{definition}
We say that a principal orbi-bundle $\mathcal P\in H^1_{orb}(\mathcal Y,\Oo_{\mathcal Y}(G))$ has rank one if there is an element $a\in \Lambda$ such that 
\begin{equation}\label{Cherncl}c^\ZZ(\mathcal P)=c_{\mathrm{fr}}\otimes a+c_{\mathrm{tors}} \end{equation}
with $0\neq c_{\mathrm{fr}}\in H_{orb}^2(\mathcal Y,\ZZ)$ being a non-torsion element, while $c_{\mathrm{tors}}$ is a torsion element. We  say additionally that a rank one principal orbi-bundle $\mathcal P$ is definite if  the image $c(\mathcal P)$ of $c^\ZZ(\mathcal P)$ via the natural map
\[H^2_{orb}(\mathcal Y,\Lambda)\rightarrow H^2_{orb}(\mathcal Y,\CC)\cong H^2(Y,\CC)\]
can be represented by a Kähler metric up to sign on $\mathcal Y$. 
\end{definition}

  Via this description, one shows, precisely like in the smooth setting (cf. \cite[Thm.~3.5]{V}, \cite[Sect.~4]{Ts99}, \cite[Prop.~4]{ist19}), that quasi-regular Vaisman manifolds correspond to definite $G$-orbi-bundles. The following is also shown in \cite[Theorem~11.9]{OV22} in a less precise form.

\begin{theorem}\label{thm: equiv qregV}
Let $G=\Lg/\Lambda$ be a one-dimensional compact complex torus. Assume that $\mathcal Y$ is a compact complex orbifold and $\pi:\mathcal P\rightarrow \mathcal Y$ is a holomorphic $G$-principal orbi-bunlde on $\mathcal Y$. Then the total space of $\mathcal P$ admits a Vaisman metric precisely when $\mathcal P$  is definite. 

In this case, there exists a finite orbifold cover $q:\mathcal Y'\rightarrow \mathcal Y$ and $\mathcal L\in\Pic^{orb}(\mathcal Y')$ such that $q^*\mathcal P\cong (\mathcal L-0_{\mathcal L})/\langle\lambda\rangle$ for some $\lambda\in \CC^*$ with $|\lambda|\neq 1$, acting by multiplication on $\mathcal L$.
\end{theorem}
\begin{proof}
We only sketch the proof, since it is analogous to the smooth case. Suppose first that the total space of the $G$ orbi-bundle $\pi:\mathcal P\rightarrow \mathcal Y$ admits a Vaisman metric $(\Omega,\theta)$ with Lee vector field $B$.
This corresponds to a collection of $\Gamma$-invariant objects $(\Omega_{\tilde U},\theta_{\tilde U},B_{\tilde U})$ on $P_{\tilde U}$ for each orbifold chart $\tilde U\in\mathcal U$ of $\mathcal Y$. Let us fix such a chart. The metric is $G$-invariant, $\Lg$ is generated by $B_{\tilde U}$ over $\CC$ and the proof of  \cite[Prop.~4]{ist19} applies to show that the Lee form determines a  $\Gamma$-invariant connection form $\alpha_{\tilde U}\in\mathcal A^1_{P_{\tilde U}}\otimes\Lg$ in the principal bundle $P_{\tilde U}$ whose curvature has the form $d\alpha_{\tilde U}=dJ\theta_{\tilde U}\otimes JB_{\tilde U}$. Since $-dJ\theta_{\tilde U}=\pi^*\omega_{\tilde U}$ for a unique Kähler form $\omega_{\tilde U}$ on $\tilde U$, the collection $(\omega_{\tilde U})_{\tilde U\in\mathcal U}$ corresponds to a global Kähler metric $\omega$ on $\mathcal Y$ such that  
\[c(\mathcal P)=[-\omega\otimes JB]\in H^2(\mathcal Y,\Lg)=H^2(\mathcal Y,\CC).\] Hence $\mathcal P$ is definite.


Note in particular that $-JB$ corresponds to an integral element of $\Lambda$ which is thus a positive multiple of a uniquely determined non-divisible element $a\in\Lambda$. This determines a canonical isomorphism $\Lg\cong \CC$ sending $a$ to $1$, and a unique element $\lambda\in\CC^*$, $|\lambda|<1$, so that we have a natural isomorphism $G\cong\CC^*/\langle\lambda\rangle$.

Let us now assume that $\mathcal P$ is definite. In this case, we can choose  $a\in\Lambda$  in \eqref{Cherncl} in a unique way so that 
$c_{\mathrm{fr}}$ is a Kähler class and $a$ is non-divisible in $\Lambda$. In particular, we can complete it to a basis $a,b$ of $\Lambda$. With respect to the isomorphism $\Lg\cong \CC$ sending $a$ to $1$, we can furthermore assume that $b\in\CC$ satisfies $\Im b>0$. Consider the diagram:
\begin{equation*}
    \begin{tikzcd}
        0\rar &\ZZ\arrow[r]\arrow[d,hook] &\CC\arrow[d,equal]\rar{\varphi} &\CC^* \rar \arrow[d, "\psi"] &0\\
        0\rar &\Lambda\rar &\CC\rar &G\rar &0
    \end{tikzcd}
\end{equation*}
where  $\varphi(z)=e^{2\pi i z}$ and $\psi$ is the induced natural projection
\[\CC^*\rightarrow \CC^*/\langle e^{2\pi i b}\rangle=G.\]
It induces the diagram on cohomology:
\begin{equation*}
    \begin{tikzcd}
       H_{orb}^1(\mathcal Y,\Oo_{\mathcal Y}^*) \arrow[d, "\psi_*"]\arrow[r, "c_1"] &H^2_{orb}(\mathcal Y,\ZZ)\dar\\
        H_{orb}^1(\mathcal Y,\Oo_{\mathcal Y}(G)) \arrow[r, "c^{\ZZ}"] & H^2_{orb}(\mathcal Y,\Lambda)
    \end{tikzcd}
\end{equation*}
where $\psi_*$ sends a line orbi-bundle $\mathcal L$ to the isomorphism class of 
\[\mathcal P':=(\mathcal L-0_{\mathcal L})/\langle e^{2\pi i b}\rangle.\]
We claim that a principal orbi-bunlde $\mathcal P'$ belongs to the image of $\psi_*$ precisely when 
\[c^\ZZ(\mathcal P')\in H^2_{orb}(\mathcal Y,\ZZ)\subset H^2_{orb}(\mathcal Y,\Lambda).\]
One direction is obvious. For the other direction, if $c^\ZZ(\mathcal P')=c_1(\mathcal L)$, it follows that $\mathcal P'=\psi_*(\mathcal L)+\mathcal P_0$, where $\mathcal P_0$ is a flat $G$-orbi-bundle 
\[ \mathcal P_0:=\tilde{\mathcal Y}\times_{\rho} G\]
with $\rho:\pi^{orb}_1(\mathcal Y)\rightarrow \Aut^{hol}_{gp}(G)$ a group homomorphism taking values in the group of holomorphic group automorphisms of $G$. But then we can lift $\rho$ to $\hat\rho:\pi^{orb}_1(\mathcal Y)\rightarrow \CC^*$ so that  $\mathcal P_0=\psi_*\mathcal L_0$ with $\mathcal L_0:=\tilde{\mathcal Y}\times_{\hat\rho}\CC$. We infer that $\mathcal P'=\psi_*(\mathcal L\otimes\mathcal L_0)$.

Let us take a finite orbifold cover $q:\mathcal Y'\rightarrow \mathcal Y$ such that $H^{orb}_1(\mathcal Y',\ZZ)$, or equivalently, $H_{orb}^2(\mathcal Y',\ZZ)$, has no torsion. Then $\mathcal P':=q^*\mathcal P\rightarrow \mathcal P$ is a finite cover and satisfies $c^\ZZ(q^*\mathcal P)=q^*c_{\mathrm{fr}}\otimes a$, hence is induced by a line orbi-bunlde $\mathcal L$ on $\mathcal Y'$ as in the above discussion. Since $\mathcal L$ is ample, it has a Hermitian structure $h$ with positive Chern curvature. The usual arguments \cite[Thm.~3.5]{V} applied over each orbifold chart then show that $\Omega:=f^{-1}i\del\ov\del f$ is a Vaisman metric on $\mathcal P'$ with Lee form $\theta:=-d\ln f$, where $f$ is the smooth real function on $\mathcal L-0_{\mathcal L}$ defined as $f(s)=||s||^2_h$.

Finally, if $F\subset\Aut(\mathcal Y')$ is the deck group of the covering $q$, then we can average $\theta$ to an $F$-invariant form $\theta_F:=\frac{1}{|F|}\sum_{g\in F}g^*\theta$. It follows that the $F$-invariant metric 
\[\Omega_F:=\frac{1}{|F|}\sum_{g\in F} g^*\Omega=-dJ\theta_F+\theta_F\wedge J\theta_F\]
is LCK. Since the action of $F$ commutes with the action of $G$ on $\mathcal P'$,  $\Omega_F$ is also $G$-invariant and the Lee vector field $B$ of $\Omega_F$ is the Lee vector field of $\Omega$. In particular, $B$ is Killing for $\Omega_F$, therefore $\Omega_F$ is a Vaisman metric which descends to  a Vaisman metric on $\mathcal P$. This concludes the proof.
\end{proof}

\begin{remark} As the base orbifold $\mathcal Y$ of a quasi-regular Vaisman manifold is required to have an integral Kähler class, the orbifold Kodaira theorem \cite{ba} implies that $\mathcal Y$ is a projective orbifold.  
\end{remark}
\begin{remark}
Note that the above theorem works well also when the total space of $\mathcal P$ is not necessarily smooth, and so corresponds to a Vaisman orbifold. Also, the only global arguments in the above proof, so which involve the orbifold condition, are the ones dealing with orbifold cohomology, in which case one transfers from the orbifold to a corresponding classifying space. 
\end{remark}
\begin{remark}
For a smooth Vaisman manifold, the base orbifold has local uniformizing groups included in a fixed torus $G$. Moreover, the above theorem shows that, up to a finite cover, the base orbifold is locally cyclic, i.e. has cyclic uniformizing groups. 
\end{remark}


\section{Deforming compact automorphism groups}\label{sec: defn}

This section deals with the deformation theory of general compact complex manifolds. The main result is the following:

\begin{theorem}\label{thm: defn}
Let $X$ be a compact complex manifold and let $H\subset\Aut_0(X)$ be a compact complex subgroup. Then for any small deformation $X_t$ of $X$ one has $\dim_\CC\Aut(X_t)\geq\dim_\CC H$. In particular, if  $\Aut_0(X)$ is compact, then for any small deformation $X_t$ of $X$, $\Aut_0(X_t)$ is a small deformation of $\Aut_0(X)$. 
\end{theorem}
\begin{proof}
Recall that if $X=(M,J)$, the small deformations $J_t$ of a the complex structure $J=J_0$ on $M$ are encoded by maps
\[\phi_t:=-\pi_t^{1,0}|_{TX^{0,1}}\in \ce(\mathcal A_X^{0,1}\otimes TX^{1,0})\]
satisfying the Maurer-Cartan equation $\ov\del\phi_t+[\phi_t,\phi_t]=0$. 
Here, $\pi_t^{1,0}:TX\otimes\CC\rightarrow TX^{1,0}_t$ denotes the $(1,0)$-projection with respect to the complex structure $J_t$, and conversely, $TX_t^{0,1}=(\id+\phi_t)(TX^{0,1})$.

Let us put
\[A^q:=\ce(\mathcal A_X^{0,q}\otimes TX^{1,0}),\quad q\in\NN. \]

The group $H$ induces natural actions on $TX^{1,0}$ and on $\mathcal A_X^{p,q}$, hence also on $A^q$, and this action commutes with the differential $\ov\del$ acting on $A^\bullet.$ Thus $H^\bullet(X,TX)$ is naturally an  $H$-representation.  Since $H$ acts holomorphically on $X$, in the sense that the map $H\times X\rightarrow X$ is holomorphic, by \cite[p. 109]{a} it follows that the $H$-represetation $H^\bullet(X,TX)$ is holomorphic.

At the same time, since $H$ is compact and acts by biholomorphisms on $X$, we can choose and fix an $H$-invariant Hermitian metric on $X$, which then induces an $H$-invariant $L^2$-scalar product on $A^\bullet.$ In particular, we can also define the $L^2$-adjoint operator $\ov\del^*$ and $\ODA:=[\ov\del,\ov\del^*]$, which will also be $H$-equivariant. Via Hodge theory, we have an $H$-equivariant isomorphism $H^\bullet(X,TX)\cong \ker\ODA$. In particular, inducing the $L^2$-scalar product on cohomology, $H^\bullet(X,TX)$ becomes a unitary $H$-representation.

Since $H$ is abelian, we can diagonalise this representation, so that it is determined by characters $\lambda:H\rightarrow \CC^*$. By the above considerations, the characters are both holomorphic and unitary functions, therefore constant. We conclude that the above representation is trivial, and in particular 
\begin{equation}\label{eq: HinvHarm}
\ker\ODA\subset (A^\bullet)^H.
\end{equation}

Recall that by \cite{kur}, there exists a germ of a complex analytic space $(\Def(X),0)$, parametrising all small deformations of $J$, with $T_0\Def(X)=H^1(X,TX)\cong\ker\ODA$. 
An explicit description of $\Def(X)$ can be given as 
\begin{equation}\label{eq: Kur}
\Def(X):=\{\phi\in A^1|\ov\del\phi+[\phi,\phi]=0,\quad \ov\del^*\phi=0\}\subset A^1\end{equation}
with the induced analytic structure form $A^1.$

\begin{claim}\label{claim}
$\Def(X)\subset (A^1)^H$.
\end{claim}
\begin{proof}
We follow the description of $\Def(X)$ given in \cite[Sect.~4]{dou}, where one can check the details. One considers appropriate Banach completions $\ov{A^\bullet}=\im\ov\del\oplus\ker\ov\del^*$ and the map $\Theta:\ov{A^1}\rightarrow\ov{A^2}$, $\phi\mapsto\ov\del\phi+[\phi,\phi]$, and defines the Banach subvariety
\[\Sigma:=\Theta^{-1}(\ker\ov\del^*)\subset\ov{A^1}\]
which is smooth at $0$ by the implicit functions theorem. Since $\Sigma$ and $\ker\ov\del^*$ meet transversely at $0$, it follows that 
$V:=\Sigma\cap\ker\ov\del^*$ is also smooth at $0$ and finite dimensional, with $T_0V=\ker\ODA|_{A^1}$. Since $V$ consists of solutions of an elliptic equation, one actually has $V\subset A^1$ and finally $\Def(X)=V\cap\Theta^{-1}(0)\subset V.$

On the other hand, the map $\Theta|_{A^1}$ is $H$-equivariant, hence one can also consider its restriction 
\[\Theta_H=\Theta:\ov{(A^1)^H}\rightarrow \ov{(A^2)^H}=(\im\ov\del)^H\oplus(\ker\ov\del^*)^H\]
together with $\Sigma^H:=(\Theta_H)^{-1}((\ker\ov\del^*)^H)$ and $V^H=\Sigma^H\cap\ker\ov\del^*\subset(A^1)^H$. By the same reasons as above, both $\Sigma^H$ and $V^H$ are smooth at $0$ and via \eqref{eq: HinvHarm} we find
\[T_0V^H=\ker\ODA|_{(A^1)^H}=\ker\ODA|_{A^1}=T_0V. \]
Since $V^H\subset V$ and both are smooth at $0$, it follows that they are equal as germs, hence $\Def(X)\subset V^H\subset (A^1)^H.$

\end{proof}

Claim \ref{claim} implies that the real Lie group $H$  acts by biholomorphisms on $(M,J')$, for any small deformation $J'$ encoded by an element $\phi\in\Def(X)$. Indeed, since $\phi$ is $H$-invariant and $TM^{0,1}_{J'}=(\id+\phi)(TX^{1,0})$, it follows that the natural action of $H$ on $TM\otimes \CC$ preserves the splitting $TM\otimes\CC=TM^{1,0}_{J'}\oplus TM^{0,1}_{J'}$, thus any element $g\in H$ satisfies $g_*J'=J'g_*$.

Let $\mathfrak{h}$ denote the Lie algebra of $H$, identified with a Lie subalgebra of $\ce(TM)$, and let us view $\mathfrak{h}^c=\mathfrak{h}+J'\mathfrak{h}$ as a complex vector space, with the complex structure induced by $J'$.  The above considerations imply that we have a natural injective $\CC$-linear morphism 
\[\mathfrak{h}^c\rightarrow H^0(TM_{J'}^{1,0}).\]
From this, we infer
\begin{align*}
\dim_\CC\Aut(M,J')&=\dim_\CC H^0(TM_{J'}^{1,0})\geq\dim_{\CC}(\mathfrak{h}+J'\mathfrak{h})=\\
&=\frac{1}{2}(2\dim_{\RR}\mathfrak{h}-\dim_{\RR}(\mathfrak{h}\cap J'\mathfrak{h}))\geq\dim_\CC H.
\end{align*}

Finally, if $H=\Aut_0(X)$, then by upper-semicontinuity  we also have 
\[\dim_{\CC}\Aut(M,J')\leq\dim_\CC H.\] This means that all inequalities in the above are equalities, so in particular $\mathfrak{h}=J'\mathfrak{h}$. Hence $J'$ endows the real Lie group $H$  with a new structure of a complex Lie group, so that $(H,J')=\Aut_0(M,J')$.  This concludes the proof.

\end{proof}

In particular, we have the immediate corollary:
\begin{corollary}\label{cor: inducedDef}
Assume that $X$ is a compact complex manifold with compact group $\Aut_0(X)$, such that $X$ is a principal $\Aut_0(X)$-orbi-bundle over an orbifold $Z$. Then any family of deformations $\mathcal X\rightarrow\Delta$ of $X$ induces family of deformations $\mathcal Z\rightarrow \Delta$ of $Z$ and $\mathcal G\rightarrow \Delta$ of $\Aut_0(X)$, such that for each $t\in\Delta$, $\mathcal X_t$ is a principal $\mathcal G_t$-orbi-bundle over $\mathcal Z_t$.
\end{corollary}

Specialising to the case of Vaisman manifolds, we obtain that the class of Vaisman manifolds with small automorphism group is open under deformations:

\begin{corollary}\label{cor: defnVaisman}
Assume $X$ is a compact complex manifold of Vaisman type with $\dim_\CC \Aut(X)=1$. Then any small deformation of $X$ is again of Vaisman type and quasi-regular.
\end{corollary}
\begin{proof}
Let $X_t$ be a deformation of $X=X_0$. By \cite[Thm.~2.6]{ov10}, $X_t$ admits an LCK metric for $t$ close to $0$. By Theorem \ref{thm: defn}, $\Aut(X_t)$ is a compact complex torus, hence \cite[Thm.~1]{ist19} implies that $X_t$ is also of Vaisman type. Finally, Theorem \ref{thm: defn} implies that $X_t$ is also quasi-regular.
\end{proof}

\section{Vaisman manifolds with $c_1(X)=0$}\label{sec: c1=0}

In this section, we start by studying the weaker Calabi-Yau condition 
\[c_1(X)=0\in H^2(X,\RR)\] 
on Vaisman manifolds. Let us recall the definition of the real Bott-Chern cohomology group
\[ H^{1,1}_{BC}(X,\RR):=\frac{\{\alpha\in\mathcal A^{1,1}_X, d\alpha=0, \alpha=\ov\alpha\}}{\{i\del\ov\del f| f\in \ce(X,\RR)\}}\]
and the fact that the identity induces a natural map 
\[l:H^{1,1}_{BC}(X,\RR)\rightarrow H^2(X,\RR)\] 
sending $c_1^{BC}(X)$ to $c_1(X)$. Thus, $c_1(X)=0$ implies that $c_1^{BC}(X)\in \ker l$.

The description of the Bott-Chern cohomology given in  \cite{IO22} immediately implies that $\ker l$ is one-dimensional, generated by the class of a transverse Kähler metric, and more particularly we have an exact sequence:
\begin{equation}\label{eq: c1}
\begin{tikzcd}
H^1(X,\RR)\rar{\psi} &H^{1,1}_{BC}(X,\RR)\rar{l} &H^{2}(X,\RR)
\end{tikzcd}
\end{equation}
where $\psi[\theta]=[-dJ\theta]$. Moreover, since the transverse Kähler classes in $H^{1,1}_{BC}(X,\RR)$ correspond to the image of the set $\mathrm{Lee}(X)$ via $\psi$, and since by \cite[Thm.~5.1]{Ts94}, $-\mathrm{Lee}(X)\cap\mathrm{Lee}(X)=\emptyset$, this allows us to split the weak Calabi-Yau manifolds into three disjoint classes:

\begin{definition}
Let $X$ be a compact Vaisman manifold with $c_1(X)=0$. We will say that:
\begin{itemize}
\item $X$ is of positive type, denoted $c_1^{BC}(X)>0$,  if $c_1^{BC}(X)$ is the class of a transverse Kähler metric, i.e. $c_1^{BC}(X)\in\psi(\mathrm{Lee}(X))$;
\item $X$ is Calabi-Yau if $c_1^{BC}(X)=0$;
\item $X$ is of negative type, denoted $c_1^{BC}(X)<0$, if $-c_1^{BC}(X)$ is the class of a transverse Kähler metric, i.e. $-c_1^{BC}(X)\in\psi(\mathrm{Lee}(X))$.
\end{itemize}
\end{definition}


%
%
%
%
%
%

For a general quasi-regular Vaisman manifold $X$ with $\mathcal Y=X/G$ and projection map $\pi:X\rightarrow\mathcal Y$, we have, via \eqref{canonical}: 
\[K_X\cong\pi^*K_{\mathcal Y}^{orb}\cong\pi^*\Oo(\phi^*(K_Y+D))\]
where $D=\sum_\alpha\left(1-\frac{1}{m_\alpha}\right)D_\alpha$ is the branch divisor of $\mathcal Y$. In particular, the following relation holds both in rational but also in Bott-Chern cohomology:
\begin{equation}\label{eq: eqBC}
c_1(X)=\pi^*(c_1(Y)-\sum_{\alpha}\left(1-\frac{1}{m_\alpha}\right)[D_\alpha]).
\end{equation}
Note that if $K_Y+D\equiv 0$ then $K_X$ is holomorphically torsion.

\begin{definition}
A compact complex  orbifold $\mathcal Y=(Y,D)$ is a Fano, Calabi-Yau, respectively of general type orbifold if 
\[c_1^{orb}(\mathcal Y)=-c_1^{orb}(K^{orb}_{\mathcal Y})=c_1(Y)-[D]\in H^2(Y,\RR)\]
is a positive, 0, respectively a negative Kähler class.
\end{definition}

  Since $\pi^*$ is injective on $(1,1)$-Bott-Chern cohomology \cite[Thm.~4.2]{IO22}, Eq. \eqref{eq: eqBC} implies that in the correspondance given by Theorem \ref{thm: equiv qregV}, the three classes also correspond to each other. In fact as we will see soon, this also holds for general Vaisman manifolds. 

Recall that any Hermitian metric $\Omega$ on a complex manifold $X$ induces a Hermitian structure on $K_X$, whose Chern curvature  can be expressed locally as 
\[\Theta_{K_X,\Omega}=\del\ov\del\log\det \Omega.\]
We also name $\rho_\Omega:=-i\Theta_{K_X,\Omega}$ \textit{the Chern-Ricci form} of $\Omega$, and it satisfies:
\[c_1^{BC}(X)=\left[-\frac{i}{2\pi}\Theta_{K_X,\Omega}\right]=\frac{1}{2\pi}[\rho_\Omega]\in H^{1,1}_{BC}(X,\RR). \]
\begin{remark}
The Calabi-Yau condition is equivalent to the existence of a function $f\in\ce(X,\RR)$ so that $\rho_\Omega=i\del\ov\del f$. In particular, as  noticed in \cite{To}, on a Calabi-Yau manifold, any Hermitian metric is conformal to a Chern-Ricci flat metric. 
\end{remark}

Recall that a Vaisman metric $(\Omega,\theta)$ on $X$ has an associated transverse Kähler metric  $\omega_0:=-dJ\theta$. Let us denote by $\rho_0$ its Ricci curvature, which can be defined as a global basic $(1,1)$-form on $X$.

\begin{lemma}\label{lem: RicciEquiv}
We have $\rho_\Omega=\rho_{0}$. In particular, $X$ is Calabi-Yau precisely when $X/\mathcal F$ is Calabi-Yau, respectively of postive or negative type precisely when $X/\mathcal F$ is Fano or of general type.
\end{lemma}
\begin{proof}
The statement is local. Pick a local transversal $Y$ for $\mathcal F$ and let $\pi:X|_Y\rightarrow Y$ be the corresponding projection, so that $\omega_0|_Y=\pi^*\eta$ for a Kähler metric $\eta$ on $Y$. If $\sigma_0$ a local frame of $K_Y$, then $\sigma:=\theta^{1,0}|_Y\wedge\pi^*\sigma_0$ is a local frame for $K_X$ over $X|_Y$. Using that $||\theta||$ is constant, we find 
\begin{gather*}
\rho_{\eta}=i\del\ov\del \log||\sigma_0||_{\eta}^2\\
 \rho_{\Omega}=i\del\ov\del \log||\sigma||^2_{\Omega}=i\del\ov\del(\log\pi^*||\sigma_0||^2_{\eta}+\log\frac{||\theta||^2}{2})=\pi^*\rho_{\eta}=\rho_{0}.\end{gather*}
\end{proof}

The above lemma tells us that all three classes of Vaisman manifolds appear in nature, simply by taking definite elliptic bundles over corresponding projective manifolds of the same sign. It also tells us that we should expect different things from Vaisman manifolds with $c_1(X)=0$, depending on the sign. 

Let us now note that the Einstein-Weyl condition implies $c_1^{BC}>0$. 
Recall first that any LCK manifold $(X,\Omega,\theta)$ has an associated minimal Kähler cover $(\hat X,\omega)$. Namely, $p:\hat X\rightarrow X$ is the minimal cover of $X$ so that $p^*\theta=d\varphi$ becomes exact, the Kähler metric is defined by $\omega:=\e^{-\varphi}p^*\Omega$, and up to constants, it only depends on the conformal class of $\Omega$.

\begin{definition}
An LCK manifold $(X,\Omega,\theta)$ is called Einstein-Weyl if its Kähler cover $(\hat X,\omega)$ is Ricci-flat. 
\end{definition}

This condition can be equivalently defined in terms of the Weyl connection, and by Gauduchon's result \cite{Gau}, LCK Einstein-Weyl manifolds are automatically Vaisman, i.e. admit a Vaisman metric in the respective conformal class. 

The Chern-Ricci form  $\rho$ of the LCK metric $g$ relates to the Chern-Ricci form $\hat\rho$ of the corresponding Kähler metric $\hat g$ by
\[\hat\rho=-i\del\ov\del \log(\det\pi^*g-n\varphi)=\pi^*(\rho-\frac{n}{2}\omega_0)\]
where $n=\dim_\CC X$. In particular, the Einstein-Weyl condition implies
\[c_1^{BC}(X)=\frac{1}{2\pi}[\rho]=\frac{n}{4\pi}[\omega_0]>0.\]

Let us also note the following:
\begin{proposition}\label{prop: KodPos}
A compact Vaisman manifold $X$ with $c_1^{BC}(X)>0$ has negative Kodaira dimension. 
\end{proposition}
\begin{proof}
Assume the contrary, so that there is some $m\geq 1$ and $0\neq \sigma\in H^0(X,K_X^{\otimes m})$. Let $D=\mathcal Z(\sigma)$ be its associated effective divisor. By hypothesis, the line bundle $K^{\otimes m}_X$ admits some Hermitian structure $h$, whose Chern curvature satisfies that
\[-\gamma:=-\frac{i}{2\pi}\Theta_{K^{\otimes m}_X,h}\]
is a  transverse Kähler metric on $X$. Let $\Omega$ be any Vaisman metric on $X$, and let $n:=\dim_\CC X$. Then \cite[eq. (1)]{Gau81} gives:
\[0<\int_D\Omega^{n-1}=\int_X\gamma\wedge \Omega^{n-1}<0\]
which is a contradiction. Thus no such section $\sigma$ exists.
\end{proof}

\section{Vaisman manifolds with $c_1^{BC}=0$}\label{sec: CY}

In this section, we collect all the main geometric properties of compact Vaisman Calabi-Yau manifolds. Most of these properties are deduced from the existence of \textit{canonical metrics}. The other main ingredient concerns their deformation theory, which was discussed in a previous section. 

\subsection{Canonical metrics}

We start by recalling the following result of Ornea and Verbitsky, which is a Yau theorem for Vaisman metrics, based on the transverse Yau theorem of El Kacimi-Alaoui \cite{ka}. We give here an equivalent reformulation of \cite[Thm~4.1]{ov22}:

\begin{theorem}\label{thm: YauV}
Let $X$ be a compact Vaisman manifold, let $G\subset\Aut(X)$ be the Lee group, let $\tau\in H^1(X,\RR)$ be a Lee class and let $\eta$ be any $G$-invariant volume form. Then there exists a unique Vaisman metric on $X$ with Lee class $\tau$ and volume form $\eta$.
\end{theorem}

Specializing to Calabi-Yau manifolds, i.e. those satisfying $c_1^{BC}(X)=0$, we obtain:

\begin{corollary}\label{thm: Ricciflat}
Let $X$ be a compact Vaisman Calabi-Yau manifold and let $\tau\in H^1(X)$ be a Lee class. Then there exists a unique Chern-Ricci flat Vaisman metric with Lee class $\tau$, up to multiplication by a positive constant.
\end{corollary}

\begin{proof}
By \cite[Proof of Thm.~5.1]{Ts94}, there exits a Vaisman metric $\Omega_0$ with Lee form $\theta_0 \in \tau$. Since $c_1^{BC}(X)=0$, $\rho_0=\mathrm{i}\partial\overline{\partial}f$ for a globally defined smooth function $f$ on $X$, where $\rho_0$ is the Chern-Ricci form of $\Omega_0$.
Then a straightforward computation gives the vanishing of $\rho_1$, the Chern-Ricci form of $\Omega_1:=e^{\frac{f}{n}}\Omega_0$, where $n=\dim_\CC X$. 

Note that since $\Omega_0$ is $G$-invariant, so is $\rho_0$. In particular, for any $\phi\in G$ we have
\[0=\phi^*\rho_0-\rho_0=i\del\ov\del(\phi^*f-f)\]
implying $\phi^*f=f$, i.e. $f$ is  $G$-invariant. Consequently, also $dv_1:=\frac{\Omega_1^{n}}{n!}$ is $G$-invariant. By Thm.~\ref{thm: YauV}, there exists a unique Vaisman metric $(\Omega,\theta)$ such that $dv_\Omega:=\frac{\Omega^n}{n!}=dv_1$ and  $\theta \in \tau$.

Since for any local trivialising section $\sigma$ of $K_X$ and any Hermitian metric $g$ one has locally
\begin{equation}\label{eq: Ricci-vol}
\rho_{g}=i\del\ov\del \log||\sigma||^2_g\quad dv_{g}=i^l||\sigma||^{-2}_g\sigma\wedge\ov\sigma\end{equation}
for some $l\in\ZZ$, it follows that the volume form determines the Chern-Ricci form of the metric. Therefore, our condition $dv_{\Omega}=dv_1$ implies $\rho_\Omega=\rho_1=0$. 

Finally, if $\Omega'$ is a different Chern-Ricci flat Vaisman metric with Lee class $\tau$, then $dv_\Omega=e^hdv_{\Omega'}$ for a global smooth function $h$. But this implies, via \eqref{eq: Ricci-vol}
\[ 0=\rho_{\Omega'}-\rho_{\Omega}=i\partial\overline{\partial}h\]
whence $h$ is constant. Therefore, due to the uniqueness part of Theorem \ref{thm: YauV}, $\Omega$ is unique up to multiplications by a positive constant.
\end{proof}

\subsection{Geometric properties}

Based on the above, we now infer a structure result for Vaisman Calabi-Yau manifolds, which can be seen as an analogue of the Beauville-Bogomolov decomposition theorem. Before stating it, let us say that for a pure orbifold $Y$ endowed with a Riemannian metric $g$, we define its holonomy as being the holonomy of the smooth locus $(Y_{\mathrm{reg}},g|_{Y_\mathrm{reg}})$.

\begin{theorem}\label{thm: holoVF}

Let $X$ be a compact Vaisman Calabi-Yau manifold. Then $\Aut_0(X)$ is a compact complex  one-dimensional Lie group and $X$ is quasi-regular. Moreover, if $g$ is a  Chern-Ricci flat Vaisman metric on $X$ and $\pi:X\rightarrow \mathcal Y=X/G$ is the associated Calabi-Yau orbifold, let  $g_0$ be the induced Ricci-flat Kähler metric on $\mathcal Y$. Then there exists a finite orbifold cover inducing a finite unramified cover of Vaisman manifolds
\begin{equation*}
    \begin{tikzcd}
    X'\arrow{r}\arrow{d}{\pi'} &X\arrow{d}{\pi}\\
    \mathcal Y'\arrow{r}{p} &\mathcal Y
    \end{tikzcd}
\end{equation*}
such that $\mathcal Y'=(Y',D=\emptyset)$ is a pure orbifold and $(Y',p^*g_0)$ splits metrically and holomorphically as
\[Y'\cong A\times\prod_{i=1}^rY_i\times\prod_{j=1}^s Z_j\]
where $A$ is a flat abelian variety, the $Y_i$'s are simply connected projective pure \Ka\ orbifolds with holonomy $\mathrm{Sp}(m_i)$ and the $Z_j$'s are simply connected projective pure \Ka\ orbifolds with holonomy $\mathrm{SU}(n_j)$.
\end{theorem}
\begin{proof}
Let us fix a Chern-Ricci flat Vaisman metric $(g,\Omega,\theta)$ on $X$ with Lee vector field $B$. Any automorphism $\phi\in\Aut_0(X)$ is homotopic to the identity, and hence it fixes the de Rham class of the Lee form $\theta$. Moreover, 
the Lee group $G$ is central in $\Aut_0(X)$ \cite[Cor.~4.5]{Ts94}, thus the LCK metric $\phi^*g$ is $G$-invariant. The Lee vector field of $\phi^*g$ is also $B$:
\[\phi^*\theta=\phi^*(g(B,\cdot))=g(B,\phi_*)=\phi^*\theta\]
hence it is Killing and $\phi^*g$ is Vaisman. It follows that $\phi^*g$ is a Chern-Ricci flat Vaisman metric with Lee class $[\theta]$ and with $\int_Xdv_{\phi^*\Omega}=\int_Xdv_\Omega$. By Theorem \ref{thm: Ricciflat}, we have thus $\phi^*g=g$, i.e. $\phi\in\Aut(X,g)$. This shows that $\Aut_0(X)=\Aut_0(X,g)$ is a compact complex Lie group, thus a torus.  Now \cite[Cor.~1]{ist19} implies that $\Aut_0(X)=G$.

The Riemannian orbifold version of the Cheeger-Gromoll theorem \cite[Thm.~2]{BZ} implies that $(\mathcal Y,g_0)$ has a finite orbifold cover $p:\mathcal Y'\rightarrow \mathcal Y$ such that 
\[(\mathcal Y', g':=p^*g_0)\cong (A,g_1)\times (\mathcal Q,g_2),\] 
where $(A, g_1)$ is a flat torus, while $(\mathcal Q, g_2)$ is a compact Ricci-flat Riemannian orbifold with $\pi_1^{orb}(\mathcal Q)=0$. We claim that since $(\mathcal Y', g')$ is Kähler, both metrics $g_1$, $g_2$ are also Kähler. For this, we need to see that if $J$ is the complex structure on $\mathcal Y'$, then $J(TA)=TA$. If this is not the case, then there exists a parallel vector field $V$ on $A$ such that $JV=W_A+W_Q$, with $W_A\in\Gamma(TA)$, $0\neq W_Q\in\Gamma(T\mathcal Q)$.
Now $W_Q$ is identified with a collection $(W_{\tilde U}\in\Gamma(T\tilde U))_{\tilde U\in\mathcal U_Q}$. Since $\nabla^{g'}(JV)=0$, it follows that also $\nabla^{g_2}W_{\tilde U}=0$ for each $\tilde U\in\mathcal U_Q$. But then the metric duals $\alpha_{\tilde U}\in\mathcal A_1(\tilde U)$ also satisfy $\nabla^{g'}\alpha_{\tilde U}=0$, and thus give rise to a non-zero harmonic one-form on the simply connected orbifold $\mathcal Q$, which is a contradiction. It follows that $J$ preserves the splitting $T\mathcal Y'=TA\oplus T\mathcal Q'$ and so induces parallel, and hence Kähler complex structures on each factor.

We take $X':=p^*X$, which is a finite cover of $X$. Since $\mathcal Q$ is simply connected, $H^2_{orb}(\mathcal Q,\ZZ)$ is torsionless, so by Theorem \ref{thm: equiv qregV} $\mathcal P:=X'|_{\{0\}\times\mathcal Q}$ is a holomorphic principal $G$-orbi-bundle of class $c(\mathcal P)=\kappa\otimes a$, where $a$ is a non-divisible element of the lattice of $G$ and $\kappa$ is a Kähler class on $\mathcal Q$. In particular, there exists $\mathcal L\in\Pic^{orb}(\mathcal Q)$ such that the holomorphic $\CC^*$-principal orbi-bundle $\mathcal Z:=\mathcal L-0_{\mathcal L}$  is a  $\ZZ$-cover of $\mathcal P$  over $\mathcal Q$. Let us further write $\kappa=d\kappa_0$, with $d\in\NN$ and $\kappa_0\in H_{orb}^2(\mathcal Q,\ZZ)$ a non-divisible class. Then there exists a degree $d$ cyclic unramified cover $\mathcal Z'\rightarrow \mathcal Z$ over $\mathcal Q$ with $c_1^{orb}(\mathcal Z')=\kappa_0$.
We have the exact sequence
\[ \begin{tikzcd}
H_2^{orb}(\mathcal Q,\ZZ)\arrow[r, "\int_{\bullet}\kappa_0"] &\ZZ=H_1(\CC^*,\ZZ)\rar &H_1^{orb}(\mathcal Z',\ZZ)\rar &H_1^{orb}(\mathcal Q,\ZZ)=0
\end{tikzcd}
\]
implying $H^{orb}_1(\mathcal Z',\ZZ)=0$, hence $H_1(Z',\ZZ)=0$. Via \cite[Thm.~4.7.5]{BG}, $\mathcal Z'\rightarrow \mathcal Q$ corresponds to a $\CC^*$-Seifert bundle with smooth total space, and so \cite[Prop.~10.2]{ko05} applies to infer  that $\mathcal Q$ is a pure orbifold.

We can thus apply the Beauville-Bogomolov decomposition theorem for pure orbifolds \cite[Thm.~6.4]{camp}  to further split $\mathcal Q$ into holonomy-irreducible factors, which thus have to be of $\mathrm{Sp}(k)$ or $\mathrm{SU}(k)$-type. Finally, as $\mathcal Y$ was projective, each factor is again projective. \end{proof}

In particular, we have:

\begin{corollary}\label{cor: KXtorsion}
If $X$ is a compact complex Calabi-Yau manifold of Vaisman type, then $K_X$ is holomorphically torsion.
\end{corollary}

Furthermore, the proof of the above theorem shows:

\begin{corollary}
Let $X$ be a compact Vaisman Calabi-Yau manifold with Lee group $G$. Then there exists a finite unramified cover $X'\rightarrow X$ and a lift of $G$ to $X'$ such that $\Sigma^G(X'):=\{x\in X'\mid G_x\neq\{1\}\}$ has complex codimension at least $2$ in $X'$.
\end{corollary}

\begin{remark}\label{rem: pure orb}
We should note that the above corollary does not hold without the assumption that $X$ is Calabi-Yau, even for $b_1(X)=1$. The  Hopf surface $X=\left(\CC^2\setminus\{0\}\right)/\ZZ$, where $\ZZ$ acts via $(z_1,z_2)\mapsto (2z_1,4z_2)$, is Vaisman quasi-regular, has $G:=\CC^*/\langle 2\rangle$ and the corresponding orbifold quotient $X/G$ is the simply connected weighted projective space $\PP(1,2)$. However,  $\Sigma^G(X')=\{z_1=0\}/\pi_1(X')$ is an elliptic curve over the orbifold point $[0:1]\in\PP(1,2)$ for any finite cover $X'$ of $X$.
\end{remark}


Concerning the deformation space of Calabi-Yau orbifolds, we have:

\begin{theorem}
    Let $X$ be a compact Vaisman Calabi-Yau manifold. Then the Kuranishi space $\Def(X)$ is a smooth universal deformation space, and any small deformation of $X$ is again Vaisman Calabi-Yau.
\end{theorem}
\begin{proof}
By Corollary \ref{cor: KXtorsion}, there exists a finite cover $p:Y\rightarrow X$ of deck group $F$, such that $K_Y$ is holomorphically trivial. Since $Y$ is also Vaisman, its Fölicher spectral sequence degenerates at $E_1$ \cite{Ts94}, therefore, by \cite[p. 154]{KKP}, \cite[Thm. 3.3]{ACRT},  the Kuranishi space $\Def(Y)$ is smooth. 

We can fix a Hermitian metric on $X$, pull it back to an $F$-invariant metric on $Y$, and define $\Def(X)$ and $\Def(Y)$ with respect to these metrics as in \eqref{eq: Kur}. Then one can identify $\Def(X)$ with the  subspace of $F$-invariant elements $\Def(Y)^F\subset\Def(Y)$. Geometrically, given a family over the disk $\pi_X:\mathcal X\rightarrow \Delta$ of central fiber $X$, one has $\pi_1(\mathcal X)\cong\pi_1(X)$ hence there is an $F$-cover $q:\mathcal Y\rightarrow \mathcal X$. The family $\pi_Y=\pi_X\circ q:\mathcal Y\rightarrow\Delta$ is then the corresponding deformation of $Y$. 

Now, as in Claim \ref{claim}, one shows that $\Def(X)\cong \Def(Y)^F\subset\Def(Y)$ is a smooth germ. Indeed, with the notations of Claim \ref{claim}, one has that $\Def(Y)=V\cap\Theta^{-1}(0)$, where
\[V=\{\phi\in A^1\mid \ov\del^*\phi=0, \quad \ov\del^*(\ov\del\phi+[\phi,\phi])=0 \}, \quad \Theta(\phi)=\ov\del\phi+[\phi,\phi].\]
But since $\Def(Y)$ is smooth at $0$ and $T_0\Def(Y)=T_0V$, we have that $(V,0)\subset(\Theta^{-1}(0),0)$ as germs. Hence also $(V^F,0)\subset(\Theta^{-1}(0),0)$, and we obtain: 
\[\Def(X)\cong\Def(Y)^F=V^F\cap\Theta^{-1}(0)=V^F\]
thus the germ $\Def(X)$ is also smooth at $0$.

Consider again a family $\pi_X:\mathcal X\rightarrow \Delta$ and the induced family $\pi_Y:\mathcal Y\rightarrow \Delta$. By  Theorem \ref{thm: defn}, for $t$ close to $0\in\Delta$ one has $H^0(X_t,TX_t)\cong H^0(X,TX).$ This implies, via \cite{Wa}, that the Kuranishi space $\Def(X)$ is universal. On the other hand, as the Frölicher spectral sequence of each $\mathcal Y_t$ degenerates at $E_1$,  $h^{n,0}(\mathcal Y_t)$ is constant in $t$. This implies that $\pi_{Y*}\Omega^n_{\mathcal Y/\Delta}$ is locally free, hence a trivial line bundle on $\Delta$. Thus one can extend any trivialisation $\sigma\in H^0(Y,K_Y)$ to a section $\sigma_{\mathcal Y}$ of $\pi_{Y*}\Omega^n_{\mathcal Y/\Delta}$, and since the condition of having constant rank is open, one can ensure, after shrinking $\Delta$, that $\sigma_{\mathcal Y}$ trivialises each $K_{\mathcal Y_t}.$ Hence each $\mathcal Y_t$, and also $\mathcal X_t$, are Vaisman Calabi-Yau manifolds.
\end{proof}

\subsection{Cohomological invariants}\label{sec: inv}

We obtain the following analytic invariants for Vaisman Calabi-Yau manifolds:
\begin{corollary}
A compact Vaisman Calabi-Yau manifold has Kodaira dimension $0$.
\end{corollary}
\begin{proof}
By Corollary \ref{cor: KXtorsion}, the canonical bundle of such a manifold is torsion, which implies the claim. 
\end{proof}

\begin{corollary}
A compact Vaisman Calabi-Yau manifold $X$ of complex dimension $n$ has algebraic dimension $n-1$.
\end{corollary}
\begin{proof}
By Theorem \ref{thm: holoVF}, $X$ fibers in elliptic curves over a projective orbifold, which implies that its algebraic dimension $a(X)\geq n-1$. At the same time, since $X$ does not satisfy the $\del\ov\del$-lemma \cite{V'}, we have $a(X)<n$ by  \cite[Cor. 5.23]{DGMS}, which implies the conclusion.
\end{proof}

Another corollary of our structure result for Vaisman Calabi-Yau manifolds is a classification of the possible fundamental groups up to finite covers, which can be seen as a topological obstruction to the existence of complex structures of Vaisman Calabi-Yau type.
\begin{corollary}
  The fundamental group of a compact Vaisman Calabi-Yau manifold contains a finite index normal subgroup which is either abelian or  $\ZZ\times H_\ZZ$, where $H_\ZZ$ is the integral Heisenberg group. Moreover, the second situation occurs only for finite quotients of the Kodaira manifolds from Example \ref{ex: Kodaira}. 
 \end{corollary}
 \begin{proof}
  Let $X$ be a compact Vaisman Calabi-Yau manifold. By Theorem \ref{thm: holoVF}, up to passing to a finite cover, we can assume that $X/G=\mathcal Y=A\times\mathcal Q$, with $A=\CC^m/\Gamma$ a compact complex torus and $\pi_1^{orb}(\mathcal Q)=0$. Then we have $\mathcal Y=A$ precisely for  Kodaira manifolds, in which case indeed $\pi_1(X)$ contains with finite index a normal subgroup of the form $\ZZ\times H_\ZZ$, cf. Example \ref{ex: Kodaira}.
  
  Let us thus assume that $\mathcal Q\neq\emptyset$. Writing as before $G=\CC/\Lambda$, we have an exact sequence:
  \[\begin{tikzcd}
  \pi_2^{orb}(\mathcal Y)\rar{\int c^\ZZ(X)}&\Lambda\rar &\pi_1(X)\rar &\pi_1^{orb}(\mathcal Y)=\Gamma\rar &0.
  \end{tikzcd}
  \]
Since $\pi_2^{orb}(\mathcal Y)=\pi_2^{orb}(\mathcal Q)=H^{orb}_2(\mathcal Q,\ZZ)$ via the Hurewicz theorem and since $c^\ZZ(X)|_{\mathcal Q}$ has rank $1$, the first map in the sequence has rank $1$. Thus,  passing to a finite cover if necessary, we can assume that we have an exact sequence
\[ \begin{tikzcd} 0\rar &\ZZ\rar &\pi_1(X)\rar &\Gamma\rar &0.\end{tikzcd}\]
In particular, this implies that $[\pi_1(X),\pi_1(X)]<\ZZ $. Since $H_1(X,\ZZ)= \pi_1(X)^{\mathrm{ab}}$ has odd rank however, we must have $[\pi_1(X),\pi_1(X)]=0$, i.e. $\pi_1(X)\cong\ZZ\times \Gamma$ is abelian.

 \end{proof}
 
 Another topological obstruction is given by the first Betti number:
 \begin{corollary}
 Let $X$ be a compact Vaisman Calabi-Yau manifold. Then
 \[ b_1(X)\leq 2\dim_\CC X-1\]
 with equality if and only if $X$ is a Kodaira manifold.
 \end{corollary}
 \begin{proof}
 Let $\Gamma\subset \pi_1(X)$ be a finite index normal subgroup which is abelian or $\ZZ\times H_\ZZ$, according to the above corollary, and let $X_\Gamma\rightarrow X$ be the corresponding cover of deck group $F$. We can also assume that $X_\Gamma/G=A\times \mathcal Q$ as above, where $G$ is the Lee group both of $X$ and of $X_\Gamma$. Then we have
 \[ b_1(X)\leq b_1(X_\Gamma)\leq 1+ 2\dim_\CC A\leq 2\dim_\CC X-1\]
 and equality can only occur if $X_\Gamma/G=A$, i.e. $X_\Gamma$ is a Kodaira manifold. Moreover, in this case $F$ acts on $\Gamma=\ZZ\times H_\ZZ$ with $\rk(\Gamma^{ab})^F=\rk(\Gamma^{ab})=2\dim_\CC X-1$. This implies that the induced action of $F$ on $\Gamma^{ab}=\ZZ\times H_\ZZ^{ab}$ is trivial. 
 
Let $\pi:X_\Gamma\rightarrow A$ be the projection map. Since the actions of $G$ and of $F$ commute, $F$ induces a natural holomorphic action on $A$, which moreover is trivial on $\pi_1(A)$. This implies that $F$ acts on $A$ by translations, so that we have a group homomorphism $\eta:F\rightarrow A$.
Assume there exists $\id\neq\gamma\in\ker\eta$. Then $\gamma$ induces a holomorphic map $f_\gamma:A\rightarrow G$ such that  $\gamma(x)=f_\gamma(\pi(x)).x$ for any $x\in X_\Gamma$, where the dot denotes the action of $G$ on $X_\Gamma$. If $f_\gamma$ is not constant, then, since $G$ is one-dimensional, $f_\gamma$ must be surjective. But  if $e$ is the identity element of $G$, then for any $x\in f_\gamma^{-1}(e)$, any $y\in\pi^{-1}(x)$ is fixed by $\gamma$, which contradicts the fact that $F$ should act freely on $X_\Gamma$. Hence $f_\gamma$ is constant and $\gamma\in G$, but this contradicts that $G$ acts effectively on $X$. Thus $\eta$ is injective, $F$ is identified with a subgroup of $A$ and $X$ is the corresponding principal $G$-bundle over the torus $A/F$, hence again a Kodaira manifold. This concludes the proof.

 \end{proof}

\section{Vaisman manifolds with $c_1^{BC}<0$}\label{sec: negative}

We turn now to the negative case. The situation is quite similar to the Calabi-Yau one. However, in order to show the existence of canonical metrics, we would need a foliated version of the Aubin-Yau theorem, giving us unique transverse negative Kähler-Einstein metrics. This probably holds, but since we were unable to find a reference, we will first show, via a Weitzenböck-type formula, that Vaisman manifolds with $c_1^{BC}<0$ are quasi-regular.

Let $X=(M,J)$ be endowed with an LCK metric $(g,\Omega,\theta)$. Recall that the Weyl connection $\mathcal{D}$ is the unique torsion free affine connection characterised by the properties
\[\mathcal DJ=0, \quad\quad \mathcal D\Omega=\theta\otimes \Omega.\]
On the universal cover of $X$, the Weyl connection lifts precisely to the Chern, or also the Levi-Civita connection of the corresponding \Ka\ metric, and we have the following relation between $\mathcal D$ and the Levi-Civita connection $\nabla^g$:
\begin{equation} \label{eq: WLC}
\mathcal D=\nabla^g-\frac{1}{2}\left(\theta\otimes\id+\id\otimes\theta-g\otimes B\right). \end{equation}

Also, $\mathcal D$ induces a natural connection on $TX^{1,0}$, which relates to the Chern connection $\nabla^{Ch}$  via (see for instance \cite[Sect.~2.3]{thesis}):
\begin{equation}\label{eq: ChernConn}
\nabla^{Ch}=\mathcal D+\theta^{1,0}\otimes\id. \end{equation}

We first show the following general formula:
\begin{lemma}
Let $X=(M,J)$ be an $n$-dimensional complex manifold with an LCK metric $(g,\Omega,\theta)$   and let $\omega_0:=-dJ\theta$. Then we have the following Weitzenb\"ock type formula:
\begin{equation}\label{eq: WeitzLCK}
    2(\ov\del^*\ov\del\xi+n\ov\del_{B^{0,1}}\xi)=\mathcal D^*\mathcal D\xi+n\mathcal D_B\xi+i(\iota_\xi\rho)^{\sharp_g}-i\frac{n}{2}(\iota_\xi\omega_0)^{\sharp_g}, \quad \xi\in\ce(X,T^{1,0}X).
\end{equation}
\end{lemma}
\begin{proof}
We let $p:\tilde X\rightarrow X$ be the universal cover, we fix $\varphi\in\ce(\tilde X,\RR)$ with $p^*\theta=d\varphi$ and we let $\tilde g:=\e^{-\varphi}p^*g$ be the corresponding \Ka\ cover. If $\mathcal D$ denotes the Levi-Civita connection on $T\tilde X$ corresponding to $\tilde g$, descending to the Weyl connection of $(\Omega,\theta)$ on $TX$, then we have the following Weitzenb\"ock formula \cite[Prop.~20.3]{mo}:
\begin{equation}\label{eq:WeitzK}
2\ov\del^{\tilde*}\ov\del\xi=\mathcal D^{\tilde *}\mathcal D\xi+i(\iota_\xi\tilde\rho)^{\sharp_{\tilde g}}, \ \ \xi\in \ce(\tilde X,T^{1,0}\tilde X)\end{equation}
where  $\tilde\rho=-i\del\ov\del\log\det\tilde g$ denotes the Ricci form of $\tilde g$, and for an operator $P$, $P^{\tilde *}$ denotes its $L^2$-adjoint with respect to the \Ka\ metric $\tilde g$. We will also use the notation $P^*$ for the $L^2$-adjoint operator of $P$ with respect to $p^*g$, respectively $g$. 

First of all, we show that 
\begin{equation}\label{eq:D*}
\mathcal D^{\tilde *}=\e^\varphi\cdot(\mathcal D^*+n\iota_B).
\end{equation}

Indeed, let $\alpha\otimes v\in\ce(\tilde X,T^*\tilde X\otimes T^{1,0}\tilde X)$, $w\in \ce(\tilde X,T^{1,0}\tilde X)$ be compactly supported. Then we find
\begin{align*}
\int_{\tilde X}\langle \mathcal D^{\tilde *}(\alpha\otimes v), w\rangle_{\tilde g}dv_{\tilde g}&=\int_{\tilde X}\langle \alpha\otimes v,\mathcal D w\rangle_{g}e^{-n\varphi}dv_g\\
&=\int_{\tilde X}\langle \alpha\otimes v, \mathcal D(\e^{-n\varphi}w)\rangle_gdv_g+\int_{\tilde X}\langle\alpha\otimes v,n\theta\otimes w\rangle_gdv_{\tilde g}\\
&=\int_{\tilde X}\langle \mathcal D^*(\alpha\otimes v),w\rangle_gdv_{\tilde g}+\int_{\tilde X}n\alpha(B)\langle v,w\rangle_gdv_{\tilde g}\\
&=\int_{\tilde X}\langle \e^\varphi\left(\mathcal D^*(\alpha\otimes v)+n\alpha(B)v\right),w\rangle_{\tilde g}dv_{\tilde g}
\end{align*}
which, by linearity, gives eq. \eqref{eq:D*}. In the same way one shows:
\begin{equation}\label{eq: ovdel*}
    \ov\del^{\tilde *}=\e^{\varphi}(\ov\del^*+n\iota_{B^{0,1}}).
\end{equation}

Replacing
\begin{equation}\label{eq:Ric}
\tilde\rho=-i\del\ov\del(\log\det p^*g-n\varphi)=p^*(\rho-\frac{n}{2}\omega_0)
\end{equation}
together with eq. \eqref{eq:D*}, \eqref{eq: ovdel*} in eq. \eqref{eq:WeitzK}, we find eq. \eqref{eq: WeitzLCK}. 

\end{proof}

\begin{theorem}\label{thm: qregNeg}
Let $X$ be a compact Vaisman manifold with $c_1^{BC}(X)<0$. Then $\Aut_0(X)=G$ is one-dimensional and $X$ is quasi-regular.
\end{theorem}

\begin{proof}
By hypothesis, $-c_1^{BC}(X)$ can be represented by a transverse Kähler form $\eta=h(J\cdot,\cdot)$. Applying the transverse Yau theorem \cite[Sect.~3.5.5]{ka}, there exists a transverse Kähler metric $\omega_0$ with transverse Ricci form $\rho_0=\eta$. But $\omega_0=-dJ\theta$ for a $G$-invariant Lee form $\theta$, where $G$ is the Lee group of $X$. Let us then fix $\Omega:=-dJ\theta+\theta\wedge J\theta$ as the fundamental form of a normalized Vaisman metric $g$ on $X$. Note that its Lee vector field satisfies
\[|B^{1,0}|_g^2=\theta^{1,0}(B^{1,0})=\frac{1}{2}.\]

 Let $\xi\in H^0(X,TX)$ be identified with a $(1,0)$-vector field. Write $\xi=\xi_0+fB^{1,0}$, where $f=2\theta^{1,0}(\xi)$ and $\xi_0\in\ker\theta^{1,0}$. Note that since $B^{1,0}$ commutes with $\xi$ \cite[Cor.~4.5]{Ts94}, $f$ must be $G$-invariant, so that in particular $[B,\xi_0]=0$ and $\ov\del_{B^{0,1}}\xi_0=0$.  Moreover, $\ov\del\xi=0$ implies
 \begin{equation}\label{eq: V1}
 \ov\del\xi_0=-\ov\del f\otimes B^{1,0}\end{equation}
 \begin{equation}\label{eq: V2}
 \ov\del f=2(\ov\del\theta^{1,0})(\xi)=i\iota_{\xi_0}\omega_0.\end{equation}

We want now to apply eq. \eqref{eq: WeitzLCK} to $\xi_0$. Using the relation \eqref{eq: WLC} between the Weyl connection  and the Levi-Civita connection, and the fact that $\nabla_B^g\xi_0=[ B,\xi_0]+\nabla_{\xi_0}^gB=0$, we find
\begin{equation}\label{eq: DB}
\mathcal D_{B^{1,0}}\xi_0=\mathcal D_B\xi_0-\ov\del_{B^{0,1}}\xi_0=\mathcal D_B\xi_0=-\frac{1}{2}\xi_0.\end{equation}
Moreover, we have:
\begin{equation}\label{eq:xi}
(\iota_{\xi_0}\omega_0)^{\sharp_g}=g(J\xi_0,\cdot)^{\sharp_g}=i\xi_0, \quad (\iota_{\xi_0}\rho)^{\sharp_g}=-ih(\xi_0,\cdot)^{\sharp_g}. \end{equation}

Replacing everything in \eqref{eq: WeitzLCK} and taking the $L^2$ product with $\xi_0$, we deduce:
\begin{align}\label{eq: Dxi}
\nonumber\int_X|\mathcal D\xi_0|^2_g dv_g&=\int_X\left( 2|\ov\del f\otimes B^{1,0}|^2_g+\frac{n}{2}|\xi_0|^2_g-|\xi_0|^2_h-\frac{n}{2}|\xi_0|^2_g   \right) dv_g\\
&\stackrel{\eqref{eq: V2}}{=}\int_X\left(|\xi_0|_g^2-|\xi_0|^2_h\right) dv_g.
\end{align}

At the same time, using again \eqref{eq: DB} and \eqref{eq: ChernConn}, we also find:
\begin{align}\label{eq: Chxi}
\nonumber\int_X|\nabla^{Ch}\xi_0|^2_gdv_g&=\int_X\left( |\mathcal D\xi_0|_g^2+2\Re g(\mathcal D \xi_0,\ov{\theta^{1,0}\otimes\xi_0})+\frac{1}{2}|\xi_0|^2_g \right)dv_g\\
\nonumber &=\int_X\left(|\mathcal D\xi_0|^2_g-\frac{1}{2}|\xi_0|^2_g\right)dv_g\\
&\stackrel{\eqref{eq: Dxi}}{=}\int_X\left(\frac{1}{2}|\xi_0|^2_g-|\xi_0|^2_h \right)dv_g.
\end{align}

Now, let us note that 
\[ (\nabla^{Ch}\xi_0)_0:=\nabla^{Ch}\xi_0-2\theta^{1,0}(\nabla^{Ch}\xi_0)\otimes B^{1,0}\]
is the orthogonal projection of $\nabla^{Ch}\xi_0$ on the orthogonal complement of $T^*X\otimes\CC B^{1,0}\subset T^*X\otimes TX^{1,0}$ with respect to the natural metric induced by $g$. This implies in particular:
\begin{align}\label{eq: V4}
 \int_X|(\nabla^{Ch}\xi_0)_0|^2_g dv_g&=\int_X\left(|\nabla^{Ch}\xi_0|^2_g-2|\theta^{1,0}(\nabla^{Ch}\xi_0)|^2_g\right)dv_g.
\end{align}

However, using eq. \eqref{eq: ChernConn}, \eqref{eq: WLC} and $\nabla^g\theta=0$, it is immediate to compute
\[\nabla^{Ch}\theta^{1,0}=\theta^{0,1}\otimes\theta^{1,0}-\frac{1}{2}g(\cdot,\pi^{1,0}).\]
We thus obtain 
\[\theta^{1,0}(\nabla^{Ch}\xi_0)=-(\nabla^{Ch}\theta^{1,0})(\xi_0)=-\frac{1}{2}g(\xi_0,\cdot)\]
and  via \eqref{eq: V4} and \eqref{eq: Chxi}, we infer:
\begin{equation*}
0\leq \int_X|(\nabla^{Ch}\xi_0)_0|^2_g dv_g=-\int_X |\xi_0|_h^2dv_g
\end{equation*}
i.e. $\xi_0=0$. In particular, via \eqref{eq: V1}, $f$ is constant and so $H^0(X,TX)=\CC B^{1,0}$. This concludes the proof.

\end{proof}

Now we are ready to show the existence of canonical metrics also in the negative case, via the orbifold version of the Aubin-Yau theorem. Before stating it, recall again the map
\[\psi:H^1(X,\RR)\rightarrow H_{BC}^{1,1}(X,\RR),\quad [\theta]\mapsto[-dJ\theta]\]
and the fact that $c_1^{BC}<0$ is equivalent to $-c_1^{BC}\in \psi(\mathrm{Lee}(X))$.

\begin{theorem}\label{thm: canMetricNeg}
Let $X$ be a compact Vaisman manifold and let $\tau\in\psi^{-1}(-2\pi c_1^{BC}(X))$.  Then there exists a  Vaisman metric $(\Omega,\theta)$ with $\theta\in\tau$ and $\rho_\Omega=dJ\theta$, which is unique up to a positive constant.
\end{theorem}
\begin{proof}
By Theorem \ref{thm: qregNeg}, $X$ is quasi-regular, and the corresponding Kähler orbifold $\mathcal Y=X/G$ satisfies $c^{orb}_1(\mathcal Y)<0$ via Lemma \ref{lem: RicciEquiv}. The orbifold Aubin-Yau theorem \cite[Thm~1.3]{F} gives a unique Kähler-Einstein metric $\eta$ on $\mathcal Y$ with $\rho_\eta=-\eta$. 

Let $\pi:X\rightarrow \mathcal Y$ be the projection map. Then $\omega_0:=\pi^*\eta$ is a transverse Kähler metric with corresponding Ricci curvature $\rho_0=-\omega_0$. In particular, 
\[\omega_0\in -2\pi c_1^{BC}(X)=\psi(\tau),\] 
hence by hypothesis there exists a closed real one-form $\theta\in\tau$ with $\omega_0=-dJ\theta$. Let us show that $\theta$ is necessarily $G$-invariant, where $G$ is the Lee group of $X$. Indeed, $G$ acts trivially on $H^1(X,\RR)$, so if $\phi\in G$, then there exists $u\in\ce(X,\RR)$ so that $\phi^*\theta=\theta+du$. But then, since $\omega_0$ is $G$-invariant, we obtain
\[0=\phi^*\omega_0-\omega_0=dd^cu \]
hence $u$ must be constant and $\phi^*\theta=\theta$.

The $(1,1)$-form $\Omega:=-dJ\theta+\theta\wedge J\theta$ is positive with $d\Omega=\theta\wedge\Omega$, and so corresponds to an LCK metric. Its Lee vector field is tangent to the Lee group $G$, therefore $(\Omega,\theta)$ is a normalized Vaisman metric. Moreover, by Lemma \ref{lem: RicciEquiv}, it satisfies $\rho_\Omega=\rho_0=dJ\theta$.

In order to show uniqueness, assume $(\Omega',\theta')$ is another normalized Vaisman metric with $\theta'=\theta+du$ and $\rho_{\Omega'}=dJ\theta'$. However, via Lemma \ref{lem: RicciEquiv}, we also have $\rho_{\Omega'}=\rho_0'$, where $\rho_0'$ is the Ricci form corresponding to the transverse Kähler metric $-dJ\theta'$. It follows that $-dJ\theta'$ is the pullback of a Kähler metric $\eta'$ on $\mathcal Y$ satisfying $\rho_{\eta'}=-\eta'$. Since such a metric is unique, it must be that $\eta'=\eta$, i.e. $dJ\theta=dJ\theta'=dJ\theta+dd^cu$, therefore $u$ is constant.  We infer that $\theta'=\theta$ and so $\Omega'=\Omega$. This concludes the proof.

\end{proof}
\begin{remark}
Note that clearly $\psi[\theta]\in -2\pi c_1^{BC}(X)$ is a necessary cohomological condition for the existence of Vaisman metrics with $\rho_\Omega=dJ\theta$, so that, unlike in the Calabi-Yau case, any Lee class $\tau$ does not solve our problem.
\end{remark}

\begin{remark}
The above also gives, after slight modifications, an alternative proof for Corollary \ref{thm: Ricciflat} of the existence of canonical metrics in the Calabi-Yau case. One needs to replace the orbifold version of the Aubin-Yau theorem with the foliated version of Yau's theorem.
\end{remark}

\begin{corollary}\label{cor: DefNeg}
Let $X$ be a compact Vaisman manifold with $c_1^{BC}(X)<0$. Then a small deformation $X_t$ of $X$ is still Vasiman with $c_1^{BC}(X_t)<0$.
\end{corollary}
\begin{proof}
Let $\mathcal X\rightarrow \Delta$ be such a deformation over the disk with $\mathcal X_0=X$. By Theorem \ref{thm: qregNeg}, $X$ is quasi-regular and has compact group $\Aut_0(X)$. By Corollary \ref{cor: inducedDef}, each $\mathcal X_t$ is a quasi-regular Vaisman manifold and there is an induced deformation of orbifolds $\mathcal Z\rightarrow\Delta$ such that for each $t\in\Delta$, $\mathcal Z_t$ is the corresponding orbifold of $\mathcal X_t$. Since $-c_1^{orb}(\mathcal Z_0)$ is a Kähler class, by the Kodaira stability theorem there exists a $2$-form $\omega$ on $\mathcal Z$ with restricts to a Kähler metric $\omega_t$ on each fiber $\mathcal Z_t$ and such that $[\omega_0]=-c_1(\mathcal Z_0)$. In particular, by continuity we also have $[\omega_t]=-c_1(\mathcal Z_t)$ for each $t\in\Delta$, therefore each $\mathcal Z_t$ is also of general type. This implies then that $\mathcal X_t$ satisfies $c^{BC}_1(\mathcal X_t)<0$.
\end{proof}

\begin{corollary}\label{cor: KodNeg}
An $n$-dimensional compact Vaisman manifold $X$ with $c_1^{BC}(X)<0$ has Kodaira dimension $n-1$.
\end{corollary}
\begin{proof}
Let $\pi:X\rightarrow\mathcal Y$ denote the projection to the corresponding orbifold. Since $\pi$ is submersive as a map of orbifolds and has parallelizable filbers, one has $K_X\cong\pi^*K_{\mathcal Y}^{orb}$.  Since $K^{orb}_\mathcal Y$ is ample, the sections of some of its power embed $\mathcal Y$ into a projective space \cite{ba}. In particular, $\mathcal Y$ has Kodaira dimension $n-1$. Since moreover $\pi$ induces injections $\pi^*:H^0(\mathcal Y,(K_{\mathcal Y}^{orb})^{\otimes m})\rightarrow H^0(X,K_X^{\otimes m})$ for any $m\in\NN$, it follows that $X$ has Kodaira dimension at least $n-1$. However, $X$ cannot have Kodaira dimension $n$, since then it would be a Moishezon manifold, and in particular it would satisfy the $\del\ov\del$-lemma \cite[Cor. 5.23]{DGMS}, contradiction with \cite{V'}. The conclusion follows. 

\end{proof}

 \section{Examples}\label{sec: examples}
Let us now give some explicit examples of Vaisman  manifolds with $c_1(X)=0$.

\begin{example}\label{ex: surfaces}
By the classification  of \cite{bel}, any  compact Vaisman surface $X$ satisfies $c_1(X)=0$. Indeed, there are three classes: 
\begin{enumerate}
\item the properly elliptic surfaces. These are elliptic orbi-bundles over curves of genus $g>1$. In particular, such surfaces  satisfy $c_1^{BC}(X)<0$;
\item the primary and the secondary Kodaira surfaces, see bellow Examples \ref{ex: Kodaira} and \ref{ex: secKod}. These satisfy $c_1^{BC}(X)=0$;
\item the diagonal Hopf surfaces and their finite quotients. These satisfy $c_1^{BC}(X)>0$. 
However, a general Hopf surface is not quasi-regular, and moreover can be deformed to a non-Vaisman surface. These examples show that the results of our paper fail for the case $c_1^{BC}>0$.
\end{enumerate}
\end{example}

\begin{example}\label{ex: Kodaira}
The simplest examples of compact Calabi-Yau orbifolds in our setting are given by the compact complex tori. 
Correspondingly, the simplest examples of compact Vaisman Calabi-Yau manifolds are given by the definite principal elliptic bundles over them, namely the Kodaira manifolds of \cite{GMPP}. These are the higher-dimensional analogues of the primary Kodaira surfaces, and have a  structure of nilmanifolds with left invariant complex structures. In particular, their fundamental groups are of the form $\Gamma=\ZZ\times\Gamma_0$ with $\Gamma_0$ a one-step nilpotent group with center $\mathcal{Z}(\Gamma_0)=\ZZ$ and $[\Gamma_0,\Gamma_0]\subset \mathcal{Z}(\Gamma_0)$ a finite index subgroup. Note that $\Gamma_0$ is an integral Heisenberg group when $[\Gamma_0,\Gamma_0]=\mathcal Z(\Gamma_0)$. 

In particular, Kodaira manifolds up to biholomorphisms are in a bijective correspondence with 
\[\{(A,\omega_0,\lambda)\mid \lambda\in\Delta^*,\quad (A,\omega_0) \text{  polarized abelian variety}\}\]
where $\Delta^*\subset\CC$ is the pointed disk. The correspondence is established as follows. Given $(A,\omega_0)$ a polarised abelian variety, $L$ a positive line bundle on $A$ with $c_1(L)=\omega_0$ and $\lambda\in\Delta^*$, $X=(L-0_L)/\langle\lambda\rangle$ is the corresponding Kodaira manifold. Moreover, any other line bundle $L'\in\Pic(A)$ with $c_1(L')=c_1(L)$ satisfies $L'=\tau_t^*L$, where $\tau_t\in\Aut(A)$ is the translation with an element $t\in A$. In particular, the induced biholomorphism  $\tilde\tau_t:L'-0_{L'}\rightarrow L-0_L$ covering $\tau_t$ is linear on the fibers and therefore commutes with the action of $\la\in\Delta^*$, thus inducing a biholomorphism of the corresponding Kodaira manifolds. Conversely, since  $H^2(A,\ZZ)$ has no torsion,  by Theorem \ref{thm: equiv qregV}, we can write any definite principal elliptic bundle over $A$ as $X=(L-0_L)/\langle\lambda\rangle$ for some positive line bundle $L\in\Pic(A)$ and a unique element $\lambda\in\Delta^*$. 
\end{example} 

\begin{example}\label{ex: secKod}
All the Calabi-Yau compact complex orbifolds of dimension $1$ can be realized as finite quotients of an elliptic curve \cite[Ex.~7.1]{F}. In particular, if they are non-smooth, then the underlying complex manifold is $\PP^1$. All of them can be realised as the base of a smooth Calabi-Yau Vaisman manifold, namely a (secondary) Kodaira surface.

We present the (adapted) construction of \cite[Sect.~V.5]{BHPV}. Assume $\mathcal Y=E/\langle\sigma\rangle$, where $E=\CC/\Lambda$ is an elliptic curve with marked point $0$ and $\sigma$ is an automorphism of $E$ of order $d$. Let $L\in\Pic(E)$ be any ample line bundle with $\sigma^*L=L$, for example $L=\mathcal O_E(0)$.  Then $\sigma$ induces a linear map $\tilde\sigma:L\rightarrow L$ covering $\sigma$, and $\tilde\sigma^d$ is a linear map covering the identity, and thus is given by fiberwise multiplication with $\mu\in\CC^*$. Taking $\nu\in\CC^*$ with $\nu^d=\mu^{-1}$, the linear map $\beta:=\nu\cdot\tilde\sigma $ is a biholomorphism of $L$ of order $d$. It follows that for any $\lambda\in\Delta^*$, the map $\rho:=\lambda\cdot\beta$ acts freely on $L-0_L$ and induces a diagram
\begin{equation*}
\begin{tikzcd}
X_1:=(L-0_L)/\langle \lambda^d\rangle \arrow[r]\arrow[d] & X_2:=(L-0_L)/\langle \rho\rangle\arrow[d]\\
E\arrow[r] &E/\langle\sigma\rangle. 
\end{tikzcd}
\end{equation*}
$X_1$ is a primary Kodaira surface, while $X_2$ is a secondary Kodaira surface. Either by the classification of compact complex surfaces or following Theorem \ref{thm: holoVF}, all Vaisman Calabi-Yau surfaces are Kodaira surfaces and are obtained in this manner.  

\end{example}

\begin{example}
For the following example, we refer to \cite[Sect.~{4.5}]{BG} concerning the general theory of weighted projective spaces. 

Let $\mathrm{w}=(w_0,\ldots,w_n)\in(\NN_{>0})^{n+1}$ be a weight with $\gcd(w_0,\ldots,w_n)=1$, and consider the weighted action $\CC^*(\mathrm{w})$ of $\CC^*$ on $\CC^{n+1}$ via 
\[\lambda.(z_0,\ldots, z_n)=(\lambda^{w_0}z_1,\ldots,\lambda^{w_n}z_n).\]
The corresponding weighted projective space is the orbifold 
\[\CC\PP(\mathrm{w}):=(\CC^{n+1}\backslash \{0\})/\CC^*(\mathrm{w}).\]
It is locally cyclic and simply connected.  Any weighted polynomial  
\begin{equation}\label{eq: f}
f\in H^0(\CC\PP(\mathrm{w}),\Oo(d))\end{equation} 
of degree $d$ defines a $\CC^*(\mathrm{w})$-invariant hypersurface $\hat X_f:=\mathcal Z(f)\backslash \{0\}\subset\CC^{n+1}\backslash \{0\}$. Assume that $\hat X_f$ is smooth. Any $\lambda\in\Delta^*$ defines a free action of $\ZZ:=\langle\lambda\rangle\subset\CC^*(\mathrm{w})$ on $\hat X_f$, so that the quotient $X_f:=\hat X_f/\ZZ$ is a compact complex manifold. 

Consider the hypersurface $Y_f=\hat X_f/\CC^*(\mathrm{w})\subset\CC\PP(\mathrm{w})$ with the inherited orbifold structure. Since $\hat X_f$ can be identified with the total space of the negative orbi-bundle $\Oo(-d)|_{Y_f}$ minus the zero section, it follows that the $G$-orbi-bundle $X_f\rightarrow Y_f$ is definite, where $G=\CC^*/\langle\lambda\rangle$, so that $X_f$ is a smooth Vaisman manifold. 

Via the adjunction formula \cite[Prop. 4.6.13]{BG}, $Y_f$ is a Calabi-Yau orbifold precisely when $d=\sum_{i=0}^nw_i$. Thus for such $d$, one obtains a smooth Vaisman Calabi-Yau manifold. For instance, the polynomials
\begin{gather*} 
f=z_0^2+z_1^6+z_2^6+z_3^6\in H^0(\CC\PP(3,1,1,1),\Oo(6))\\
h=z_0^2z_1+z_1^5+z_2^5+z_3^5\in H^0(\CC\PP(2,1,1,1),\Oo(5))
\end{gather*}
define smooth Vaisman Calabi-Yau threefolds. The corresponding orbifold $Y_f$ is a smooth K3 surface, while $Y_h$ is an orbifold K3 surface with orbifold singularity $[1,0,0,0]$  of uniformiszing group $\ZZ/2\ZZ$.

This procedure can be easily generalised to complete intersections and produces many examples. The list of all weighted K3 hypersurfaces was found by Reid \cite[Sect.~4]{Reid}, and more general weighted Calabi-Yau orbifolds are found in abundance for instance in \cite{IF},\cite{CLS}, which thus can be used to produce  Vaisman Calabi-Yau manifolds as above.  In \cite[Sect.~8]{IF}, conditions on $(\mathrm{w},d)$ are also given to ensure that a general polynomial $f$ as in \eqref{eq: f} defines a smooth  $\hat X_f$, and so a smooth Vaisman manifold.

\end{example}


\begin{thebibliography}{100}

\bibitem[A12]{a} D. Akhiezer, {\it Lie Group Actions in complex Analysis}, Vol. 27. Springer Science \& Business Media, 2012.

\bibitem[ACRT18]{ACRT} B. Anthes, A. Cattaneo, S. Rollenske, A. Tomassini, \textit{$\del\ov\del$-Complex symplectic and Calabi–Yau manifolds: Albanese map, deformations and period maps}, Ann. Glob. Anal. Geom. 54 (2018), 377--398.

\bibitem[Ba57]{ba} W. L. Baily, {\it On the imbedding of V-manifolds in projective spaces}, Amer. J. Math. 79 (1957), 403--430.

\bibitem[BHPV04]{BHPV} W. P. Barth, K. Hulek, C. A. M. Peters, A. Van de Ven, {\it Compact complex surfaces}, volume 4 of {\it Ergebnisse der Mathematik und ihrer Grenzgebiete. 3. Folge}. Springer-Verlag, Berlin, second edition, 2004.

 \bibitem[Be84]{Be84} A. Beauville, {\it Variétés Kähleriennes dont la première classe de Chern est nulle}, J. Differential Geom. 18, no. 4 (1984), 755--782.
 
 \bibitem[B00]{bel} F. Belgun, {\it On the metric structure of non-K\"ahler complex surfaces}, Math. Ann. 317 (2000), 1--40.


\bibitem[BG07]{BG} C. P. Boyer, K. Galicki, {\it Sasakian Geometry}, Oxford Mathematical Monographs, Oxford University Press, Oxford, 2007.

\bibitem[BZ94]{BZ} J. Borzellino, S. Zhu, {\it The Splitting Theorem for Orbifolds}, Illinois J. of Math. 38 (1994), 679--691.

\bibitem[Ca04]{camp} F. Campana, {\it Orbifoldes \`a premi\`ere classe de Chern nulle}, in \textit{The Fano Conference},  pp. 339--351. University Torino, Torino, 2004.

\bibitem[CLS90]{CLS} P. Candelas, M. Lynker, R. Schimmrigk, {\it Calabi-Yau manifolds in weighted $\PP_4$}, Nuclear Phys. B 341 no. 2, (1990) 383--402. 


\bibitem[DGMS75]{DGMS} P. Deligne, Ph. Griffiths, J. Morgan,J D. Sullivan, {\it Real homotopy theory of Kähler manifolds}, Invent. Math. 29(3) (1975), 245--274.

\bibitem[Do66]{dou} A. Douady, {\it Le probl\`eme des modules pour les vari\'et\'es analytiques complexes}, in {\it Séminaire Bourbaki}, Vol. 9 (1964/1965), no. 277, Soc. Math. France, Paris, 1995, 7--13.





\bibitem[E90]{ka}  A. El Kacimi-Alaoui,  {\it Op\' erateurs transversalement elliptiques sur un feuilletage riemannien et applications}, Compositio Math. 73 (1990), 57--106.


\bibitem[Fa18]{F} M. Faulk, {\it On Yau's theorem for effective orbifolds}, Expositiones Mathematicae, 37.4 (2019), 382--409.

\bibitem[Ga81]{Gau81} P. Gauduchon, {\it Le théorème de dualité pluriharmonique}, C.R. Acad. Sci. Paris 293 (1981), 59--63.


\bibitem[Ga95]{Gau} P. Gauduchon, {\it Structures de Weyl-Einstein, éspaces de twisteurs et variétés de type $S^1\times S^3$}, J. reine angew. Math. 469 (1995), 1--50. 

\bibitem[GMPP04]{GMPP} G. Grantcharov, C. McLaughlin, H. Pedersen, and Y.S. Poon, {\it Deformations of Kodaira manifolds}, Glasg. Math. J., 46(2) (2004), 259--281.

\bibitem[IF00]{IF} A. R. Iano-Fletcher, {\it Working with weighted complete intersections}, Explicit birational geometry of 3-folds, London Math. Soc. Lecture Note Ser., vol. 281, pp. 101--173. Cambridge Univ. Press, Cambridge, 2000. 

\bibitem[Is18]{thesis} N. Istrati, {\it Conformal structures on compact complex manifolds}, Ph.D. Thesis, Université Sorbonne Paris Cité, 2018.

\bibitem[Is19]{ist19} N. Istrati, {\it Existence criteria for special locally conformally Kähler metrics}, Annali di Matematica Pura ed Applicata (1923 -) volume 198 (2019), 335--353. 

\bibitem[IO22]{IO22} N. Istrati, A. Otiman, {\it Bott-Chern cohomology of compact Vaisman manifolds},  to appear in Trans. Amer. Math.  Soc. DOI: 10.1090/tran/8832.

\bibitem[KKP08]{KKP} L. Katzarkov, M. Kontsevich, T. Pantev, {\it Hodge theoretic aspects of mirror symmetry}. In {\it From Hodge theory to integrability and TQFT $tt^*$-geometry}, volume 78 of Proc. Sympos. Pure Math., pp. 87--174. Amer. Math. Soc., Providence, RI, 2008.



\bibitem[Ko05]{ko05} J. Koll\'ar, {\it Einstein metrics on five-dimensional Seifert bundles}, J. Geom. Anal. 15 no. 3  (2005), 445--476. 

\bibitem[Ku62]{kur} M. Kuranishi, {\it On the locally complete families of complex analytic structures}, Ann. Math. (2) 75 (1962) 536--577.

\bibitem[Mo07]{mo} A. Moroianu, {\it Lectures on Kähler Geometry}, London Mathematical Society Student Texts 69. Cambridge University Press, Cambridge, 2007. 

\bibitem[Na75]{Na} I. Nakamura, {\it Complex parallelisable manifolds and their small deformations}, J. Differential Geom. 10 (1975), 85--112.







\bibitem[OV10]{ov10} L. Ornea, M. Verbitsky, \textit{Locally conformal Kähler manifolds with potential}, Math. Ann. 348 (2010), 25--33.


\bibitem[OV22A]{ov22} L. Ornea, M. Verbitsky, {\it A Calabi-Yau theorem for Vaisman manifolds}, arXiv: 2206.08808.

\bibitem[OV22B]{OV22} L. Ornea, M. Verbitsky, \textit{Principles of Locally Conformally Kähler Geometry}, arXiv:2208.07188. 


\bibitem[Re80]{Reid} M. Reid, {\it Canonical 3-folds}, Journ\'ees de G\'eom\'etrie Alg\'ebrique d'Angers, Juillet 1979/Algebraic Geometry, Angers, 1979, pp. 273–310. Sijthoff \& Noordhoff, Alphen aan den Rijn, 1980.

\bibitem[Ro11]{Ro} S. Rollenske, {\it The Kuranishi space of complex parallelisable nilmanifolds}, J. Eur. Math. Soc. (JEMS), 13(3) (2011), 513--531.



\bibitem[To15]{To} V. Tosatti, {\it Non-Kähler Calabi-Yau manifolds}, Analysis, complex geometry, and mathematical physics: in honor of Duong H. Phong 644 (2015), 261--277.

\bibitem[TW10]{TW} V. Tosatti, B. Weinkove, {\it The complex Monge-Ampère equation on compact Hermitian manifolds}, J. Amer. Math. Soc. 23,  no.4 (2010), 1187--1195.

\bibitem[Ts94]{Ts94}  K. Tsukada, {\it Holomorphic forms and holomorphic vector fields on compact generalized Hopf manifolds}, Compositio Math. 93 (1) (1994), 1--22. 

\bibitem[Ts97]{Ts97} K. Tsukada, {\it Holomorphic Maps of Compact Generalized Hopf Manifolds}, Geometriae Dedicata 68 (1997), 61--71.

\bibitem[Ts99]{Ts99} K. Tsukada, {\it The canonical foliation of a compact generalized Hopf manifold}, Differential Geom. Appl. 11 (1999), 13--28.


\bibitem[Va80]{V'} I. Vaisman, {\it On Locally and Globally Conformal K\"ahler Manifolds}, Trans. Amer. Math. Soc. 262, no. 2 (1980), 553--542. 

\bibitem[Va82]{V} I. Vaisman, {\it Generalized Hopf manifolds}, Geom. Dedicata 13 (1982), 231--255.

\bibitem[Wa69]{Wa} J. Wavrik, {\it Obstructions to the existence of a space of moduli}. In \textit{Global Analysis (Papers in Honor of K. Kodaira)}, pp. 403--414. Univ. Tokyo Press, Tokyo, 1969.

\bibitem[Ya78]{Yau} S. T. Yau, {\it On the Ricci curvature of a compact Kähler manifold and the complex Monge-Ampère equation}, I. Comm. Pure Appl. Math 31, no.3 (1978), 339--411.








\end{thebibliography}
\end{document}